\documentclass[a4paper]{article}
\usepackage[affil-it]{authblk}
\usepackage{amssymb}
\usepackage{amsfonts}
\usepackage{amsmath}
\usepackage{multirow}
\usepackage{float}
\usepackage{xcolor}
\usepackage{bm}
\usepackage{caption}
\usepackage{subcaption}
\usepackage{array}
\usepackage{etoolbox}
\usepackage{mathtools}
\usepackage{dutchcal}
\usepackage{graphicx}
\usepackage{booktabs, makecell, multirow}
\usepackage{geometry}
\usepackage[export]{adjustbox}
\usepackage{graphbox}
\usepackage{setspace}
\usepackage{bookmark}
\usepackage{threeparttable}
\usepackage{tikz}

\newcommand{\stabregion}{
\raisebox{2pt}{\tikz{\draw[-,solid,line width = 1.0pt](0.,0.8mm) -- (5.5mm,0.8mm)}}\,\,\,}
\newcommand{\eigenwithout}{\tikz{\draw[magenta,fill] (2.75mm,0.75mm) circle (0.65mm);}\,\,\,}
\newcommand{\eigenwith}{\tikz{\draw[blue,fill] (2.75mm,0.75mm) circle (0.65mm);}\,\,\,}

\newcommand{\exactsol}{
{\raisebox{0pt}{\tikz{\draw[fill] (2.5mm,0.75mm) circle (0.65mm);\draw[black,line width = 1.0pt](0.,0.8mm) -- (5mm,0.8mm);}}}}
\newcommand{\numericalsol}{
{\raisebox{0pt}{\tikz{\draw[fill,magenta] (2.5mm,0.75mm) circle (0.65mm);\draw[magenta,line width = 1.0pt](0.,0.8mm) -- (5mm,0.8mm)}}}}

\newcommand{\hybridCFL}{
{\raisebox{0pt}{\tikz{\node[rectangle,draw,yscale=0.6,rotate=45, fill, green]{}}}}}
\newcommand{\rkdCFL}{
{\raisebox{0pt}{\tikz{\draw[fill, red] (2.5mm,0.75mm) circle (0.65mm);\draw[red,line width = 1.0pt](0.,0.8mm) -- (5mm,0.8mm);}}}}

\let\originalleft\left
\let\originalright\right

\renewcommand{\left}{\mathopen{}\mathclose\bgroup\originalleft}
\renewcommand{\right}{\aftergroup\egroup\originalright}
\renewcommand*{\arraystretch}{1.0}

\newcommand{\SR}{\textrm{S\hskip-1pt R}}
\newcommand{\rks}{\mathcal{s}}

\newcolumntype{L}[1]{>{\raggedright\let\newline\\\arraybackslash\hspace{0pt}}m{#1}}
\newcolumntype{C}[1]{>{\centering\let\newline\\\arraybackslash\hspace{0pt}}m{#1}}
\newcolumntype{R}[1]{>{\raggedleft\let\newline\\\arraybackslash\hspace{0pt}}m{#1}}

\spacing{1.0}

\newtheorem{myremark}{Remark}
\newtheorem{prop}{Property}
\newfont{\notapolice}{cmss10}
\begin{document}

%\begin{frontmatter}
\title{An \textit{a posteriori} strategy for adaptive schemes in time and space}

\author{Maria T. Malheiro%
  \thanks{Electronic address: \texttt{mtm@math.uminho.pt}}}
\affil{Centre of Mathematics and Department of Mathematics, University of Minho\\
  Campus of Azur\'em, 4800-058 Guimar\~aes, Portugal}

\author{Gaspar~J.~Machado%
  \thanks{Electronic address: \texttt{gjm@math.uminho.pt}}}
\affil{Centre of Physics and Department of Mathematics, University of Minho\\
  Campus of Azur\'em, 4800-058 Guimar\~aes, Portugal}

\author{St{\'e}phane Clain%
  \thanks{Electronic address: \texttt{clain@math.uminho.pt}}}
\affil{Centre of Physics and Department of Mathematics, University of Minho\\
  Campus of Azur\'em, 4800-058 Guimar\~aes, Portugal}

\date{October, 2020}

\maketitle

\begin{abstract}
A nonlinear adaptive procedure for optimising both the schemes in
time and space is proposed in view of increasing the numerical
efficiency and reducing the computational time. The method is based
on a four-parameter family of schemes we shall tune in function of
the physical data (velocity, diffusion), the characteristic size in
time and space, and the local regularity of the function leading to
a nonlinear procedure. The \textit{a posteriori} strategy we adopt
consists in, given the solution at time $t^n$, computing a candidate
solution with the highest accurate schemes in time and space for all
the nodes. Then, for the nodes that present some instabilities, both
the schemes in time and space are modified and adapted in order to
preserve the stability with a large time step. The updated solution
is computed with node-dependent schemes both in time and space. For
the sake of simplicity, only convection-diffusion problems are
addressed as a prototype with a two-parameters five-points finite
difference method for the spatial discretisation together with an
explicit time two-parameters four-stages Runge-Kutta method. We
prove that we manage to obtain an optimal time-step algorithm that
produces accurate numerical approximations exempt of non-physical
oscillations.
\end{abstract}
%\begin{keyword}
%optimal time step \sep stability \sep finite difference method \sep MOOD \sep high-order \sep non-stationary convection- diffusion
%\end{keyword}
%\end{frontmatter}

%\noindent{\bf MSC:}XXXXXX.

\section{Introduction}
\label{sec:introduction}

High order discrete schemes for equations involving hyperbolic
operators are likely to produce non physical oscillations. Even for
linear problems, a nonlinear routine is mandatory to control the
over- and under-shooting in the vicinity of points where the
solution presents large gradients. Technologies such as MUSCL or
WENO \cite{MUSCL,WENO}, among the most popular, have been developed
for half a century and manage to efficiently reduce or eliminate the
numerical instabilities, namely to avoid the oscillations near the
discontinuities and the extrema points. More specifically, for the
convection-diffusion problems, different approaches have been
considered: variable-order non-oscillatory scheme (VONOS),
hybrid-linear parabolic approximation (HLPA), sharp and monotonic
algorithm for realistic transport (SMART), weighted-average
coefficient ensuring boundedness (WACEB), convergent and universally
bounded interpolation scheme for the treatment of advection
(CUBISTA) and an adaptive bounded version of the QUICK with
estimated streaming terms (QUICKEST) called ADBQUICKEST (see
references for these methods in \cite{Ferreira}). The numerical
solution obtained with these schemes is at least second-order
accurate in regions where the solution is smooth enough but retains
the first-order approximation in regions where the solution presents
large gradients for the sake of stability.

All the stabilisation procedures mainly address the scheme in space
whereas the scheme in time is merely discretised with a Runge-Kuta
(RK) method  or its Strong Stability Preserving (SSP) version
\cite{Wu88,WuOs88}. Very little attention has been paid on the time
discretisation and a nonlinear dynamical procedure for optimising
both the schemes in time and space is desirable to increase the
numerical efficiency. Several traditional numerical methods were
revisited in order to articulate space and time schemes together,
aiming to increase the allowable time step.
Bourchtein~\cite{Bourchtein} constructed an explicit central
diffe\-rence method of second order applied to the one- and
two-dimensional advection equations based on the generalised
leap-frog type method with the main goal of increasing the allowable
time step, with some deterioration in the accuracy of the solution.
Chadha and Madden~\cite{Chadha} consider the numerical solution of a
linear time dependent advection–diffusion problem by an implicit
two-weight, three-point central finite difference scheme. They
extend the scheme proposed by them in \cite{Chadha 2}, to
incorporate an optimal time step selection algorithm for the method.
The resulting method, based on optimal values of weights and optimal
time-stepping, is of fifth-order in space, and third-order in time.

Other numerical schemes are based on prediction-correction
techniques  \cite{Verfurth}.  In \cite{Araya}, the authors present
an adaptive finite element scheme for the
advection-reaction-diffusion equation based on a stabilized finite
element method combined with a residual error estimator. The
adaptive process corrects the meshes allowing to capture boundary
and inner layers very sharply and without significant oscillations.
In \cite{Knopp}, the authors have developed some \textit{a
posteriori} error estimates using a stabilized scheme combined with
a shock-capturing technique to control the local oscillations in the
crosswind direction. In \cite{Wang}, the author has introduced an
error estimate for the advection–diffusion equation based on the
solution of local problems on each element of the triangulation.
Artificial compression method (ACM) based filter scheme is also
investigated in  \cite{Yee}. The spatially fourth-order (or higher
non-dissipative) scheme in space is used at all times but additional
numerical dissipation is made at the shock layers to control the
stability.

In the last decade, the \textit{a posteriori} paradigm has been
developed to prevent the numerical solution from creating non
physical oscillations: the so-called Multidimensional Optimal Order
Detection Method (MOOD) method  \cite{Clain,Clain2}. The principle
consists in building a candidate solution with the highest order
scheme in space. Then the guess solution admissibility is analysed
using detectors to check some physical properties of non physical
oscillations. The nodes that present non-compliance's solution are
tagged and the numerical scheme is only altered for that points by
reducing the scheme order (basically adding more viscosity or
reducing the polynomial reconstruction degree). The numerical
approximation for the cured nodes and their neighbours is computed
again to eliminate the oscillations. The goal of this the paper is
the design of an adaptive technique, based on the \textit{a
posteriori} paradigm, to provide the optimal choice of the time and
space schemes, that enables to compute a stable solution with the
largest time step.

To present our strategy, we deal with the one-dimension linear
convection-diffusion problem since all the important ingredients are
already addressed in this simple scalar equation and enable to
highlight the connections between the scheme in time and the scheme
in space. Notice that, even a linear problem requires a nonlinear
routine for stabilisation when dealing with rough solutions and the
convection-diffusion equation involves the two main operator in
simulations. For the sake of simplicity, we only consider two
schemes in space and two schemes in time with very different
characteristics and we aim to optimise the space accuracy and the
time step to produce physically admissible numerical solutions with
no spurious oscillations.

The crucial point lies on the confrontation between the
discretisations in space and time: the eigenvalues associated to the
space scheme have to fit into the stability domain of the RK scheme.
Several parameters play a major role in the trade-off between
accuracy and stability. On the one hand, the scheme in space
together with the  P\'eclet number determine the eigenvalues
distribution in the complex plane. On the other hand, the RK scheme
stability region and the accuracy is characterised by the Butcher
tableau \cite{Butcher}. At last, the stability condition between the
two schemes is controlled by the Courant–Friedrichs–Lewy (CFL)
condition depending on the time discretisation parameter $\Delta t$.

In this study, we consider a two-parameter family of finite
difference schemes in space that  characterise the spectrum and the
accuracy. Similarly, the time discretization is a two-parameter
family of schemes by considering a four-stage RK method where we
impose to be at least second-order. The global method is a
four-parameter family of schemes we shall tune in function of the
physical data (velocity, diffusion), the characteristic size in time
and space, and the local regularity of the function leading to a
nonlinear procedure since the scheme depends on the approximation
itself.

The design of the schemes in space has to respect some basic
principles in order to produce an eligible blending. First we only
consider five-point finite difference methods involving the same
centred stencil but the coefficients are node dependent. Similarly
the four-stage RK scheme is also node dependent and it is mandatory
that the four sub-steps are the same for all the nodes for the sake
of compatibility leading to some additional constraints in the
design of the Butcher's tableau. The \textit{a posteriori} strategy
we adopt consists in, given the solution at time $t^n$, computing a
candidate solution with the highest accurate scheme in time and
space for all the nodes. Then applying the detectors, we determine
the nodes to be corrected and switch the scheme in space.
Unfortunately, it usually leads to a reduction of the local time
step for the sake of stability. To overcome such a problem and
preserve a large time step, we may have to switch the time scheme.
We then obtain, the optimal combination of time and space schemes
for each node and iteration while we use the common time step to
update the solution.

The organisation of the paper is the following. We briefly present
in Section~2 the spectral analysis of the five-point finite
difference method focusing on the spectral curves description with
respect to the two free parameters. Section~3 is dedicated to the
four-stage RK method where we analyse the impact of the two free
parameters on the stability region. The stability of the time and
space schemes combination is carried out in Section~4 where we
determine the optimal time step for each scenario. Finally,
Section~5 presents our \textit{a posteriori} method to design a node
by node optimal scheme both in time and space. At last in Section~6
conclusions and perspectives are drawn.

\section{Discrete convection-diffusion operator analysis}
\label{sec:space}

This section is dedicated to the spectrum of the discrete
convection-diffusion operator we shall use in our stability
analysis. Let $\phi=\phi(x)$ be a smooth 1-periodic function defined
in $\mathbb R$, \textit{i.e.}, $\phi(x+1)=\phi(x)$. We define the
convection-diffusion operator
\begin{equation}
\label{ConvDiff_operator} \mathfrak E[\phi]=-u\phi'+\kappa \phi'',
\end{equation}
where $u\geq 0$ and $\kappa\geq 0$ are the convective and diffusive
coefficients, respectively. We restrict the study case to the
bounded domain $[0,1]$ by applying a periodic condition. Let
$I\in\mathbb N$ and $\Delta x=1/I$. We denote $x_i=i\Delta x$,
$i\in\mathbb Z$, a uniform discretisation of the real axis and set
$\Psi=(\psi_i)_{i\in\mathbb{Z}}$ a vector with an infinite number of
real value entries. Periodicity yields $\psi_{i+I}=\psi_i$, for all
$i\in\mathbb{Z}$, hence the relevant data is only given by
components $\psi_i$, $i=1,\ldots,I$. We shall use the same notation
$\Psi=(\psi_i)_{i=1}^{I}$ to denote both the whole vector and its
finite representation, being the other components given by
periodicity.

\subsection{Five-point discrete schemes}

A generic conservative five-point numerical scheme is defined by an
ordered list of 5 coefficients which we indicate with
\begin{equation}
\label{generic_scheme_E} E=(a_{-2},a_{-1},a_{0},a_{1},a_{2}), \quad
a_{j}\in\mathbb R, j=-2,\ldots,2,
\end{equation}
where the coefficients satisfy the null summation constraint for the
sake of conservation
\[
a_{-2}+a_{-1}+a_{0}+a_{1}+a_{2}=0.
\]
For any $I$-periodic vector $\Psi=(\psi_i)_{i\in\mathbb Z}$, the
$E$-scheme applied to $\Psi$ provides the vector $E\Psi$ given
component-wise by
\[
(E\Psi)_i=a_{-2}\psi_{i-2}+a_{-1}\psi_{i-1}+a_0\psi_{i}+a_{1}\psi_{i+1}+a_{2}\psi_{i+2},
\quad i\in \mathbb{Z}.
\]
Periodicity of $\Psi$ yields that $E\Psi$ also satisfies the
periodicity property hence the finite vector representation
$E\Psi=((E\Psi)_i)_{i=1}^{I}$ completely describes the whole vector.
We highlight the four particular cases
\begin{alignat*}{2}
&
E_1=\left(\frac{1}{12},-\frac{2}{3},0,\frac{2}{3},-\frac{1}{12}\right),
\qquad
&& E_2=\left(-\frac{1}{12},\frac{4}{3},-\frac{5}{2},\frac{4}{3},-\frac{1}{12}\right),\\
&  E_3=\left(-\frac{1}{2},1,0,-1,\frac{1}{2}\right), &&
E_4=\left(1,-4,6,-4,1\right),
\end{alignat*}
that provide approximations for the first-, second-, third-, and
fourth-order derivatives, respectively. Denoting
$\Psi=(\phi(x_i))_{i\in \mathbb Z}$ for any regular $1$-periodic
function $\phi=\phi(x)$, consistency errors read
\begin{alignat*}{2}
&  (E_1\Psi)_i=\phi^{(1)}(x_i)+\mathcal O(\Delta x^4), \qquad
&& (E_2\Psi)_i=\phi^{(2)}(x_i)+\mathcal O(\Delta x^4),\\
&  (E_3\Psi)_i=\phi^{(3)}(x_i)+\mathcal O(\Delta x^2), &&
(E_4\Psi)_i=\phi^{(4)}(x_i)+\mathcal O(\Delta x^2).
\end{alignat*}
The fourth-order five-point discretisation of operator
\eqref{ConvDiff_operator} is given by the optimal combination
\[
E=-\frac{u}{\Delta x}\left(E_1-\frac{E_2}{Pe}\right),
\]
where $Pe=\frac{u\,\Delta x}{\kappa}$ represents the cell P\'eclet
number. Instabilities may appear and one has to damp the
oscillations by using third- and fourth-derivative approximations.
To this end, we consider more general five-point conservative
schemes of the form
\begin{equation}
\label{discEq} E=E(\theta_3,\theta_4,Pe)=-\frac{u}{\Delta
x}\left(E_1-\frac{E_2}{Pe}+\theta_3E_3+\theta_4E_4\right),
\end{equation}
parameterised by  $\theta_3,\theta_4\in \mathbb R$, and $Pe$.
Substituting the expressions of $E_1$, $E_2$, $E_3$, and $E_4$,
scheme $E$ reads
\begin{multline}
\label{schemeE}
E(\theta_3,\theta_4,Pe)=-\frac{u}{\Delta x}\left(\frac{(12\theta_4-6\theta_3+1)Pe+1}{12Pe},-\frac{(12\theta_4-3\theta_3+2)Pe+4}{3Pe},\right.\\
\left.\frac{12\theta_4
Pe+5}{2Pe},-\frac{(12\theta_4+3\theta_3-2)Pe+4}{3Pe},\frac{(12\theta_4+6\theta_3-1)Pe+1}{12Pe}\right).
\end{multline}

\begin{myremark}
If $u=0$ and $\kappa\ne 0$ (pure diffusive problem), $Pe=0$ and
expression~\eqref{schemeE} should be rewritten specifically for
$u=0$ taking into account $\theta_3$ and $\theta_4$. If $\kappa=0$
and $u\ne 0$ (pure convective problem) we will write $Pe=+\infty$,
and expression~\eqref{schemeE} should be rewritten for $\kappa=0$
taking again into account $\theta_3$ and $\theta_4$.
\end{myremark}

%{\bf remark}
%We do not mention the scale coefficient $\frac{u}{\Delta x}$ at that stage but will be fundamental in the time discretisation situation by involving the $CFL$ constraint.

\subsection{Spectra}

\label{sec2-spectra} Due to the periodicity assumption,
scheme~\eqref{generic_scheme_E} results in a circulant square matrix
$A$ of order $I$ with entries
\[
a_{ij}=
\begin{cases}
a_{j-i} & \text{if } |j-i|\leq 2\\
a_{j-i-I} & \text{if } |j-i|\geq I-2\\
0  & \text{otherwise}.
\end{cases}
\]
Schematically, we have
\[
A=\begin{bmatrix}
a_0 & a_1 & a_2 & 0 & 0 & 0 & 0 & \cdots & 0 & a_{-2} & a_{-1}\\
a_{-1} & a_0 & a_1 & a_2 & 0 & 0 & 0 & \cdots & 0 & 0 & a_{-2}\\
a_{-2} & a_{-1} & a_0 & a_1 & a_2 & 0 & 0 & \cdots & 0 & 0 & 0\\
0 & a_{-2} & a_{-1} & a_0 & a_1 & a_2 & 0 & \cdots & 0 & 0 & 0\\
\vdots & \vdots & \vdots & \vdots & \vdots & \vdots & \vdots & \ddots & \vdots & \vdots & \vdots\\
a_2 & 0 & 0 & 0 & 0 & 0 & 0 & \cdots & a_{-1} & a_0 & a_1\\
a_1 & a_2 & 0 & 0 & 0 & 0 & 0 & \cdots & a_{-2} & a_{-1} & a_0\\
\end{bmatrix}.
\]
The eigenvectors of the circulant matrices associated to schemes
$E_1$, $E_2$, $E_3$, and $E_4$ are the same, hence they also are the
eigenvectors of the circulant matrix associated to scheme $E$ and
independent of $\theta_3$, $\theta_4$, and $Pe$. They are given by
\[
v^{(i)}=\begin{bmatrix} 1 & w_i & w_i^2 & \cdots & w_i^{I-1}
\end{bmatrix}, \quad i=1,\ldots,I,
\]
with $w_i=\exp\left (2\pi \mathfrak i (i\Delta x)\right )$ and
$\mathfrak i$ the unit imaginary number \cite{Gantmacher}.
Identifying the discrete scheme $E(\theta_3,\theta_4,Pe)$ with the
respective circulant matrix $A(\theta_3,\theta_4,Pe)$, the
eigenvalue $\lambda_i(E(\theta_3,\theta_4,Pe))$ associated to
$v^{(i)}$ depends on coefficients $\theta_3$, $\theta_4$, and $Pe$.
Taking into account scheme~\eqref{discEq} we get
\begin{equation}\label{eigenvalues_A}
\lambda_i(E(\theta_3,\theta_4,Pe))=
%\lambda_i\left(-\frac{u}{\Delta x}\left(E_1-\frac{E_2}{Pe}+\theta_3\,E_3+\theta_4\,E_4\right)\right)=
-\frac{u}{\Delta x}\left
(\lambda_i(E_1)-\frac{\lambda_i(E_2)}{Pe}+\theta_3\lambda_i(E_3)+\theta_4\lambda_i(E_4)\right).
\end{equation}
The eigenvalues for the four particular operators read:
\begin{itemize}
\item $\displaystyle \lambda_i(E_1)=\mathfrak{i}\sin(2\pi i\Delta x)\left(1-\frac{1}{3}\Big\{\cos(2\pi i\Delta x)-1\Big\}\right)\in \left[-\mathfrak{i}\frac{4}{3},\mathfrak{i}\frac{4}{3}\right]$,
\item $\displaystyle \lambda_i(E_2)=\big (\cos(2\pi i\Delta x)-1\big)\left(2-\frac{1}{3}\Big\{\cos(2\pi i\Delta x)-1\Big \}\right)\in\left[-\frac{16}{3},0\right]$,
\item $\displaystyle \lambda_i(E_3)=2\mathfrak{i}\sin(2\pi i\Delta x)\big(\cos(2\pi i\Delta x)-1\big)\in[-2\mathfrak{i},2\mathfrak{i}]$,
\item $\displaystyle \lambda_i(E_4)=4\big(\cos(2\pi i\Delta x)-1\big)^2\in[0,16]$.
\end{itemize}
Letting $I\rightarrow+\infty$, we obtain the continuous
parameterised spectral curve
\begin{equation}
\label{equ_lambda}
\lambda(\rks;\theta_3,\theta_4,Pe)=\frac{u}{\Delta
x}\rho(\rks;\theta_3,\theta_4,Pe), \quad \rks\in [0,1],
\end{equation}
with
\begin{equation}
\label{equ_rho}
\rho(\rks;\theta_3,\theta_4,Pe)=x(\rks;\theta_3,\theta_4,Pe)+\mathfrak{i}y(\rks;\theta_3,\theta_4,Pe)
\end{equation}
where
\begin{align}
& x(\rks;\theta_3,\theta_4,Pe)=\frac{1}{Pe}\Big\{\cos(2\pi \rks)-1\Big \} \left[2-\Big\{\cos(2\pi \rks)-1\Big\} \left(\frac{1}{3}+4Pe\theta_4\right)\right],\label{x_equation}\\
& y(\rks;\theta_3,\theta_4,Pe)=-\sin(2\pi \rks)\left [
1-\Big\{\cos(2\pi \rks)-1\Big\}
\left(\frac{1}{3}-2\theta_3\right)\right ].\label{y_equation}
\end{align}
\begin{myremark}
It is worth noting that the spectral curve shape is characterised by
function $\rho(\rks;\theta_3,\theta_4,Pe)$ while $\frac{u}{\Delta
x}$ is a scaling factor we shall blend with the time step parameter
to produce a CFL-like coefficient.
\end{myremark}

\subsection{Centered and upwind schemes}

The \textbf{centered scheme} corresponds to $\theta_3=0$ and
$\theta_4=0$ and provides the optimal fourth-order of accuracy. The
corresponding scheme for finite and positive $Pe$ reads
\[
E_\text{c}=-\frac{u}{\Delta
x}\left(\frac{1+Pe}{12\,Pe},-\frac{4+2Pe}{3\,Pe},\frac{5}{2\,Pe},-\frac{4-2Pe}{3\,Pe},\frac{1-Pe}{12\,Pe}\right).
\]
We present in Table~\ref{tab:spectra_centered} the spectral curves
for the centered scheme when $Pe\in]0,+\infty[$ and also for $Pe=0$
and $Pe=+\infty$.

\begin{table}[!ht]\centering
\caption{Spectral curves --- centered scheme: $\theta_3=0$ and
$\theta_4=0$.} \label{tab:spectra_centered}
\begin{tabular}{@{}ccc@{}}\toprule
$Pe$ & spectrum curve $\lambda$ & $\rho_\text{c}$ \\\midrule
$]0,+\infty[$ & \makecell[l]{
$\displaystyle\lambda(\rks;Pe)=\frac{u}{\Delta x}\left(x_\text{c}(\rks;Pe)+\mathfrak{i}y_\text{c}(\rks;Pe)\right)$\\[0.25cm]
$\displaystyle x_{\text{c}}(\rks;Pe)=\frac{1}{Pe}\left(\cos(2\pi \rks)-1\right)\left(2-\frac{1}{3}\Big(\cos(2\pi \rks)-1\Big)\right)$\\[0.25cm]
$\displaystyle y_\text{c}(\rks;Pe)=-\sin(2\pi
\rks)\left(1-\frac{1}{3}\Big(\cos(2\pi \rks)-1\Big)\right)$} &
\includegraphics[height=3cm,valign=m]{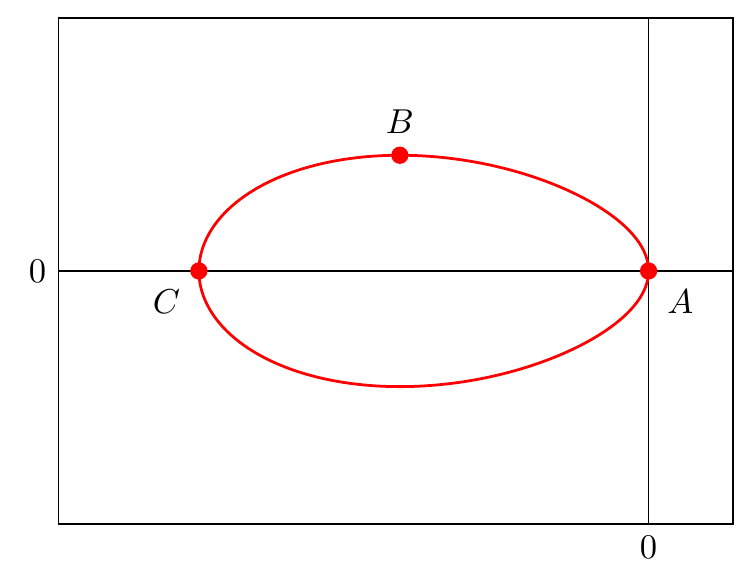} \\\midrule
$0$ & \makecell[l]{
$\displaystyle \lambda(\rks)=\frac{\kappa}{\Delta x^2}\left(x_\text{c}(\rks)+\mathfrak{i}y_\text{c}(\rks)\right)$\\[0.25cm]
$\displaystyle x_\text{c}(\rks)=(\cos(2\pi \rks)-1) \left[2-\frac{1}{3}(\cos(2\pi \rks)-1)\right]$\\
$\displaystyle y_\text{c}(\rks)=0$} &
\includegraphics[height=3cm,valign=m]{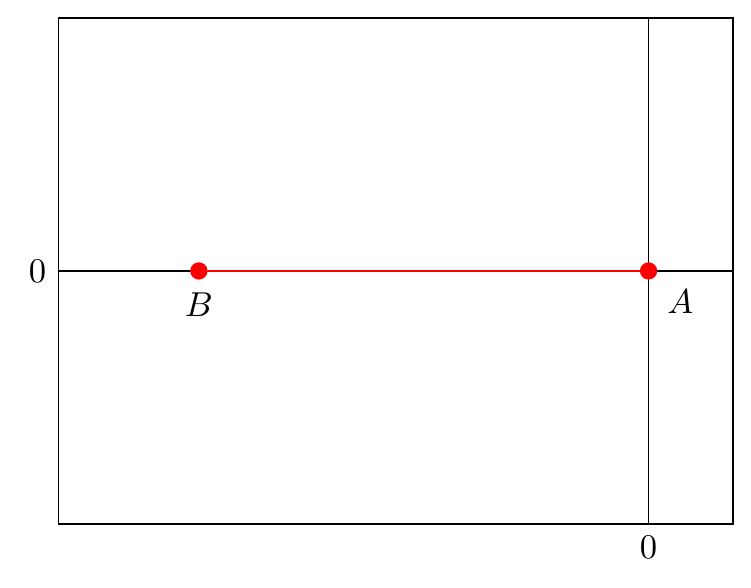} \\\midrule
$+\infty$ & \makecell[l]{
$\displaystyle \lambda(\rks)=\frac{u}{\Delta x}\left(x_\text{c}(\rks)+\mathfrak{i}y_\text{c}(\rks)\right)$\\[0.25cm]
$\displaystyle x_\text{c}(\rks)=0$\\[0.25cm]
$\displaystyle y_\text{c}(\rks)=-\sin(2\pi
\rks)\left[1-\frac{1}{3}\Big(\cos(2\pi \rks)-1\Big)\right]$} &
\includegraphics[height=3cm,valign=m]{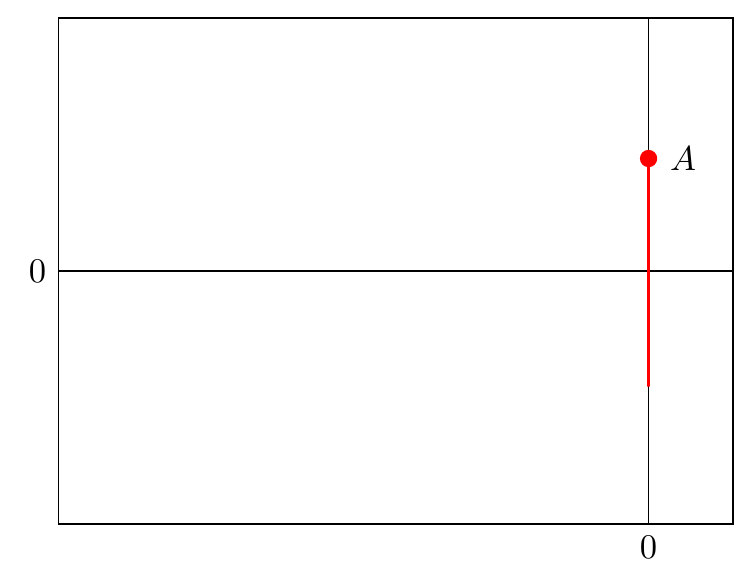} \\\bottomrule
\end{tabular}
\end{table}

We define the \textbf{weak upwind scheme} by cancelling coefficient
$a_2$ since we have assumed $u\geq 0$. Hence the following relation
has to be satisfied
\begin{equation}\label{weak_upwind}
\theta_4+\frac{\theta_3}{2}=\frac{Pe-1}{12\,Pe}.
\end{equation}
In order to get the optimal accuracy, we set $\theta_3=0$ leading to
a third-order and deduce by relation \eqref{weak_upwind} that
$\theta_4=\frac{Pe-1}{12Pe}$. The scheme for finite and positive
$Pe$ then reads
\[
E_\text{w}=-\frac{u}{\Delta
x}\left(\frac{1}{6},-\frac{Pe+1}{Pe},\frac{Pe+4}{2Pe},\frac{Pe-3}{3Pe},0\right).
\]
We present in Table~\ref{tab:spectra_weak} the spectral curves for
the weak-upwind scheme for $Pe\in]0,+\infty[$ and also for $Pe=0$
and $Pe=+\infty$.
\begin{table}[ht!]\centering
\caption{Spectral curves --- weak upwind scheme: $\theta_3=0$,
$\theta_4=\frac{Pe-1}{12Pe}$.} \label{tab:spectra_weak}
\begin{tabular}{@{}ccc@{}}\toprule
$Pe$ & spectrum curve $\lambda$ & $\rho_\text{w}$ \\\midrule
$]0,+\infty[$ & \makecell[l]{
$\displaystyle \lambda(\rks;Pe)=\frac{u}{\Delta x}\left(x_\text{w}(\rks;Pe)+\mathfrak{i}y_\text{w}(\rks;Pe)\right)$\\[0.25cm]
$\displaystyle x_\text{w}(\rks;Pe)=\frac{1}{Pe}\left(\cos(2\pi \rks)-1\right)\left(2-\frac{Pe}{3}\Big\{\cos(2\pi \rks)-1\Big \} \right)$\\[0.25cm]
$\displaystyle y_\text{w}(\rks;Pe)=-\sin(2\pi
\rks)\left(1-\frac{1}{3}\left(\cos(2\pi \rks)-1\right)\right)$} &
\includegraphics[height=3cm,valign=m]{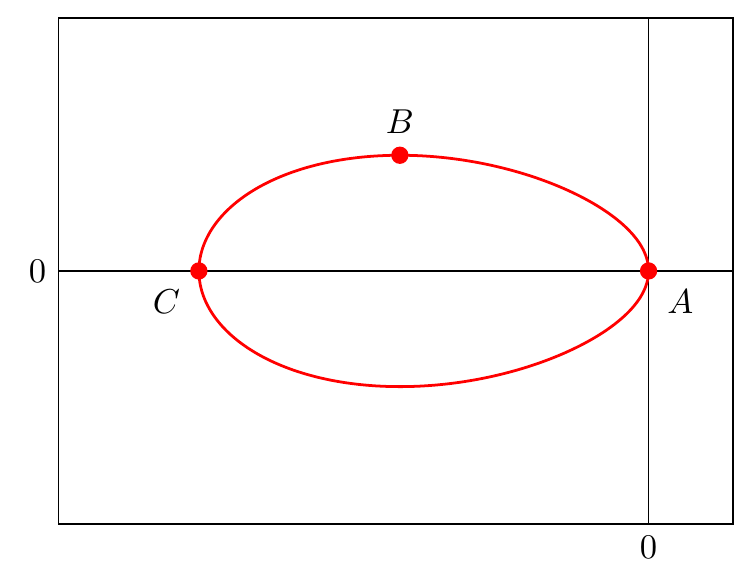} \\\midrule
$0$ & \makecell[l]{
$\displaystyle \lambda(\rks)=\frac{\kappa}{\Delta x^2}\left(x_\text{w}(\rks)+\mathfrak{i}y_\text{w}(\rks)\right)$\\[0.25cm]
$\displaystyle x_\text{w}(\rks)=2 \Big(\cos(2\pi \rks)-1\Big)$\\[0.25cm]
$\displaystyle y_\text{w}(\rks)=0$} &
\includegraphics[height=3cm,valign=m]{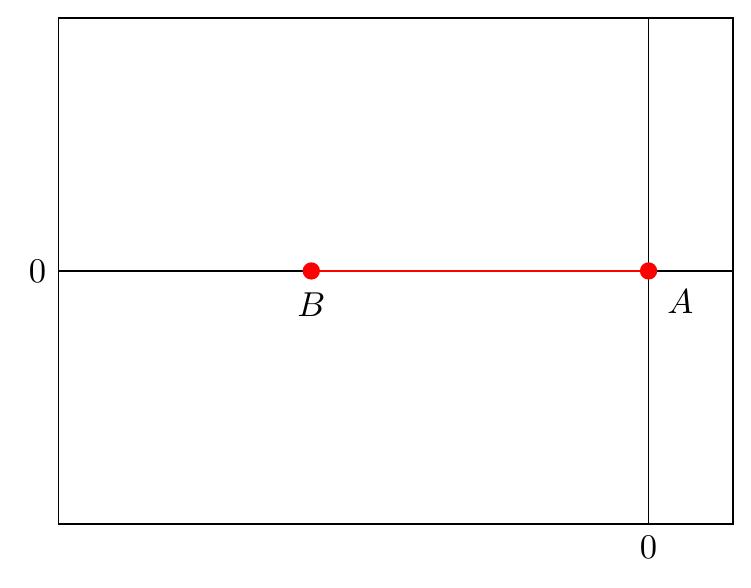} \\\midrule
$+\infty$ & \makecell[l]{ $\displaystyle \lambda(\rks)=
\frac{u}{\Delta x}\left(x_\text{w}(\rks)+\mathfrak{i}y_\text{w}(\rks)\right)$\\[0.25cm]
$\displaystyle x_\text{w}(\rks)=-\frac{1}{3}(\cos(2\pi \rks)-1)^2$\\[0.25cm]
$\displaystyle y_\text{w}(\rks)=-\sin(2\pi
\rks)\left(1-\frac{1}{3}\left(\cos(2\pi \rks)-1\right)\right)$} &
\includegraphics[height=3cm,valign=m]{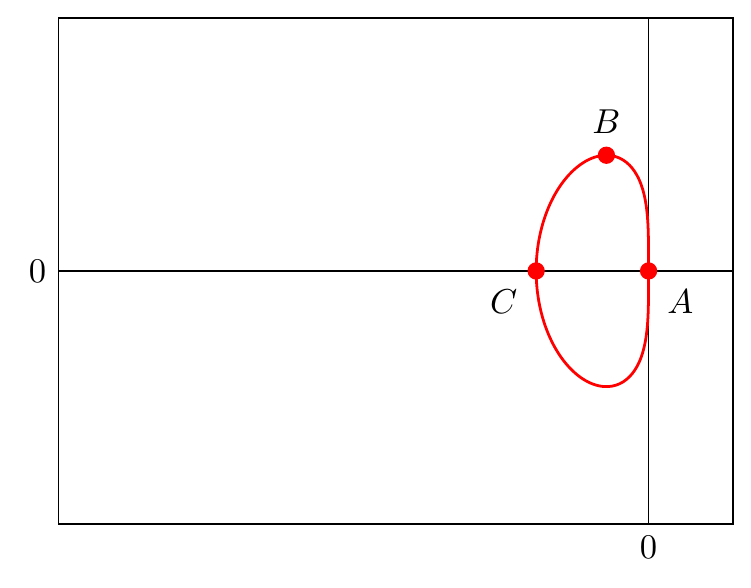} \\\bottomrule
\end{tabular}
\end{table}

The \textbf{strong upwind scheme} consists in cancelling both
coefficients $a_1$ and $a_2$ leading to the second-order full upwind
three-point scheme with
\[
\theta_3=\frac{3-Pe}{3\,Pe}\quad\text{and}\quad
\theta_4=\frac{3\,Pe-7}{12\,Pe}.
\]
The corresponding scheme for finite and positive $Pe$ reads
\[
E_\text{s}=-\frac{u}{\Delta
x}\left(\frac{Pe-2}{2Pe},-\frac{2Pe-2}{Pe},\frac{3Pe-2}{2Pe},0,0\right).
\]
We present in Table~\ref{tab:spectra_strong} the spectral curves for
the strong upwind scheme when $Pe\in]0,+\infty[$ and also for $Pe=0$
and $Pe=+\infty$.
\renewcommand{\arraystretch}{2}
\begin{table}[ht!]\centering
\caption{Spectral curves --- strong upwind scheme:
$\theta_3=\frac{3-Pe}{3Pe}$, $\theta_4=\frac{3Pe-7}{12Pe}$.}
\label{tab:spectra_strong}
\begin{tabular}{@{}ccc@{}}\toprule
$Pe$ & spectrum curve $\lambda$ &  $\rho_\text{s}$
\\\midrule
$]0,+\infty[$ & \makecell[l]{ $\displaystyle
\lambda(\rks;Pe)=\frac{u}{\Delta
x}\left(x_\text{s}(\rks;Pe)+\mathfrak{i}y_\text{s}(\rks;Pe)\right)$
\\[1cm]
$\displaystyle x_\text{s}(\rks;Pe)=\frac{1}{Pe}\left(\cos(2\pi \rks)-1\right)\left(2-(Pe-2)\Big\{\cos(2\pi \rks)-1\Big\}\right)$ \\[1cm]
$\displaystyle y_\text{s}(\rks;Pe)=-\sin(2\pi
\rks)\left(1-\frac{Pe-2}{Pe}\Big\{\cos(2\pi \rks)-1\Big\}\right)$}
& \makecell{$]0,1[$  \\
\includegraphics[height=3cm,valign=m]{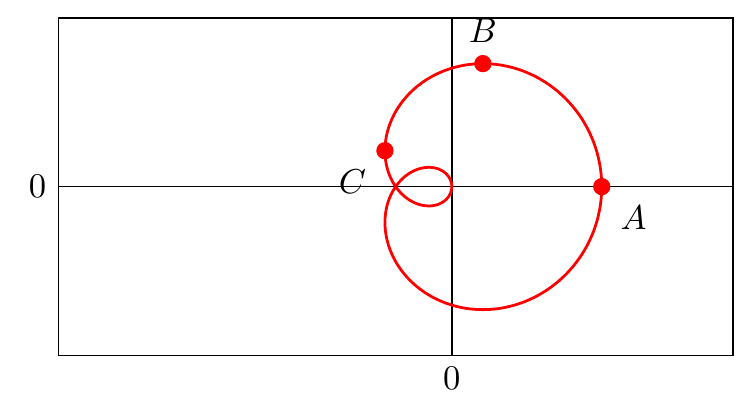}
\\ $[1,+\infty[$\\
\includegraphics[height=3cm,valign=m]{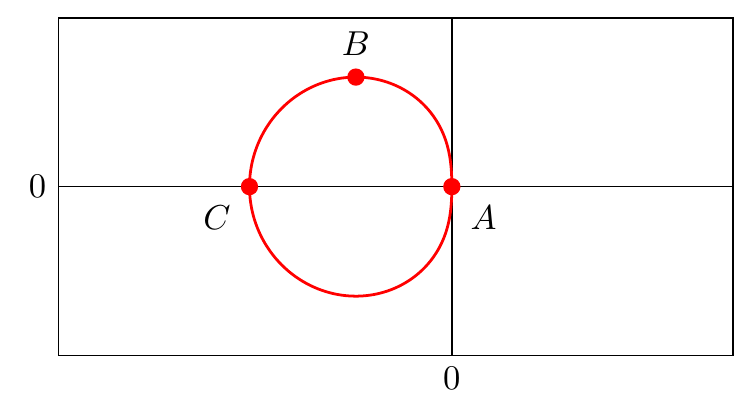}}
\\\midrule

$0$ & \makecell[l]{ $\displaystyle
\lambda(\rks)=\frac{\kappa}{\Delta
x^2}\rho(\rks)=\frac{\kappa}{\Delta x^2}\left(x(\rks)+\mathfrak{i}y(\rks)\right)$ \\[0.25cm]
$\displaystyle x_\text{s}(\rks)=2\cos(2\pi \rks)\Big(\cos(2\pi \rks)-1\Big)$ \\[0.25cm]
$\displaystyle y(\rks)=-2\sin(2\pi \rks)\Big(\cos(2\pi\rks)-1\Big)$}
&
\includegraphics[height=3cm,valign=m]{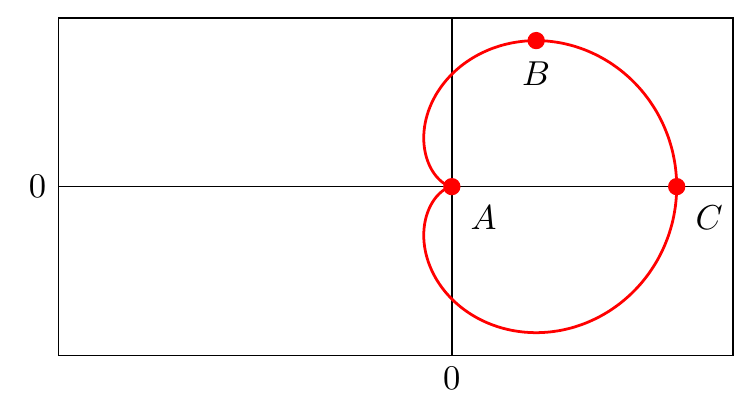}\\\midrule

$+\infty$ & \makecell[l]{$\displaystyle
\lambda(\rks)=\frac{u}{\Delta x}\rho(\rks)=
\frac{u}{\Delta x}\left(x_\text{s}(\rks)+\mathfrak{i}y_\text{s}(\rks)\right)$ \\[0.25cm]
$\displaystyle x_\text{s}(\rks)=-(\cos(2\pi \rks)-1)^2$ \\[0.25cm]
$\displaystyle y_\text{s}(\rks)=-\sin(2\pi \rks)\left [ 2-\cos(2\pi
\rks)\right]$}  &
\includegraphics[height=3cm,valign=m]{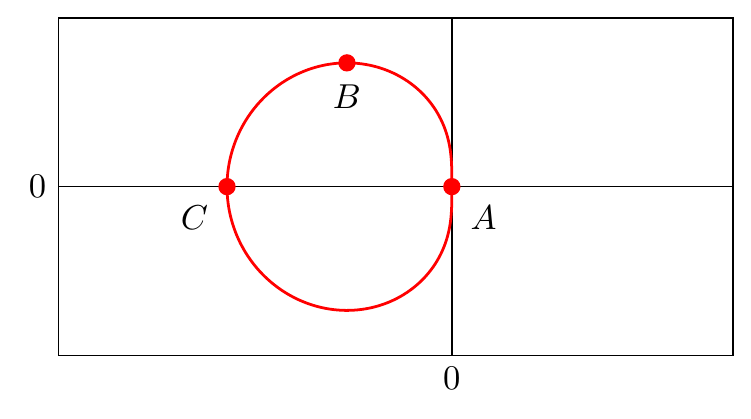}
\\\bottomrule
\end{tabular}
\end{table}
\renewcommand{\arraystretch}{1.0}

Table~\ref{tab:spectra} reports the relevant points of the spectral
curves of the centered, weak upwind and strong upwind schemes marked
on Tables~\ref{tab:spectra_centered}, \ref{tab:spectra_weak} and
\ref{tab:spectra_strong}, respectively.

\begin{table}[ht!]\centering
\begin{threeparttable}
\caption{Relevant points of the spectral curves} \label{tab:spectra}
\begin{tabular}{@{}ccccc@{}}\toprule
$Pe$ & scheme & $A$ & $B$ & $C$ \\ \midrule $]0, +\infty[$ &
\makecell{centered \\ weak upw} & $(0,0)$&\makecell{$(\frac{-2\sqrt{6}-1}{2Pe},1.37)$\\[0.25cm]$(-0.5-\frac{2.45}{Pe},1.37)$} &
\makecell{$(-\frac{16}{3Pe},0)$\\[0.25cm]$(-\frac{4(3+Pe)}{3Pe},0)$ }\\ \cmidrule{3-5}
\makecell{$]0,1[$ \\
[0.25cm]$[1,+\infty[$} & strong upw & \makecell{$(\frac{4(1-Pe)}{Pe},0)$ \\
[0.25cm]$(0,0)$} &
\makecell{$(x_{\text{s},B,1}, y_{\text{s},B,1})\tnote{a}$ \\
[0.25cm]$(x_{\text{s},B,2}, y_{\text{s},B,2})\tnote{b}$} &
\makecell{$(\frac{1}{Pe(Pe-2)},\frac{(Pe-1)\sqrt{3-2Pe}}{Pe(Pe-2)})$\\[0.25cm]$(\frac{4(1-Pe)}{Pe},0)$}

 \\\midrule
$0$ &
\makecell{centered\\[0.25cm]weak upw\\[0.25cm]strong upw} &
\makecell{$(0,0)$} &
\makecell{$(-5.33,0)$\\[0.25cm]$(-4,0)$\\[0.25cm]$(1.5,2.6)$} &
\makecell{--\\[0.25cm]--\\[0.25cm]$(4,0)$}
 \\\midrule
$+\infty$ &
\makecell{centered\\[0.25cm]weak upw\\[0.25cm]strong upw} &
\makecell{$(0,1.37)$\\[0.25cm]$(0,0)$\\[0.25cm]$(0,0)$} &
\makecell{--\\[0.25cm]$(-0.5,1.37)$\\[0.25cm]$(-1.87,2.2)$} &
\makecell{--\\[0.25cm]$(-1.33,0)$\\[0.25cm]$(-4,0)$}
 \\\bottomrule
\end{tabular}
\begin{tablenotes}
      \item[a]
      $x_{\text{s},B,1}=\frac{(3-Pe+\sqrt{\xi})(1+Pe-\sqrt{\xi})}{4Pe(Pe-2)},y_{\text{s},B,1}=\frac{\sqrt{-2(2Pe-3)-2(Pe-1)\sqrt{\xi}}(3Pe-3-\sqrt{\xi})}{4Pe(Pe-2)}$,
      \item[b] $x_{\text{s},B,2}=\frac{(3-Pe-\sqrt{\xi})(1+Pe+\sqrt{\xi})}{4Pe(Pe-2)},y_{\text{s},B,2}=
\frac{\sqrt{-2(2Pe-3)+2(Pe-1)\sqrt{\xi}}(3Pe-3+\sqrt{\xi})}{4Pe(Pe-2)}$,
\\ where $\xi=3Pe^2-10Pe+9$.
    \end{tablenotes}
\end{threeparttable}
\end{table}

We compare in Table~\ref{fig:centered_weak_upwind} the spectral
curves for the centered, weak upwind, and strong upwind schemes when
$Pe=1, 5, 30$. We observe the similarity between centered and weak
upwind space discretisations.

\begin{table}[!ht]\centering
\caption{Spectral curves for different values of $Pe$ for the three
pairs $\theta_3, \theta_4$ cases.} \label{fig:centered_weak_upwind}
\begin{tabular}{@{}cccc@{}}\toprule
scheme & $\theta_3$ & $\theta_4$ & spectra\\\midrule
\makecell{centered\\(fourth-order)} & $0$ & $0$ &
\includegraphics[valign=m]{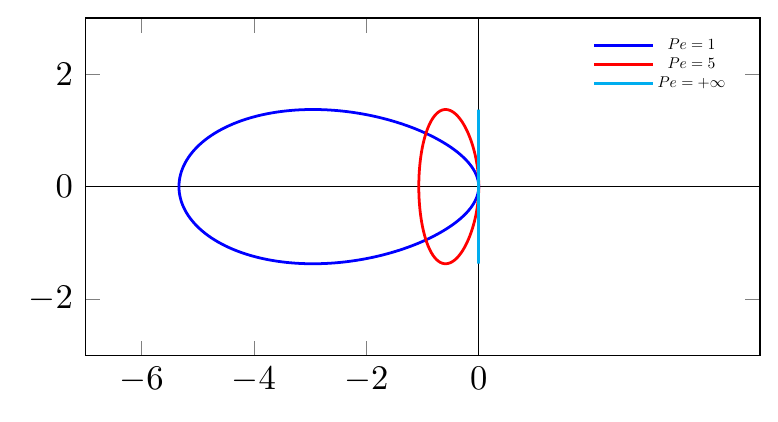}\\\hline
\makecell{weak upwind\\(third-order)} & $0$ & $\frac{Pe-1}{12Pe}$ &
\includegraphics[valign=m]{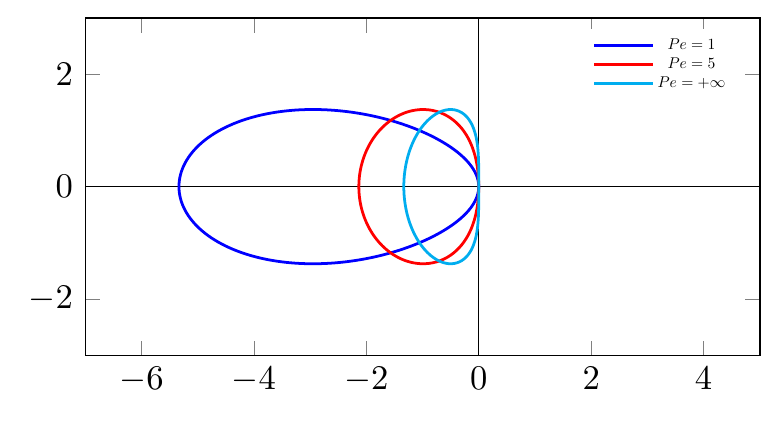}\\\hline
\makecell{strong upwind\\(second-order)} & $\frac{3-Pe}{3\,Pe}$ &
$\frac{3\,Pe-7}{12\,Pe}$ &
\includegraphics[valign=m]{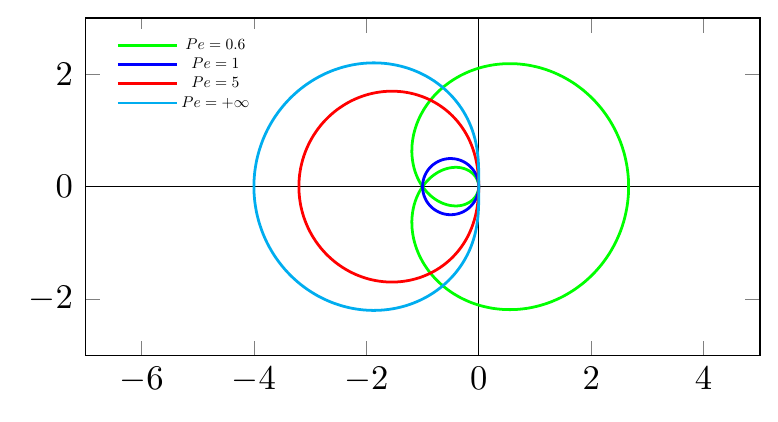}\\\bottomrule
\end{tabular}
\end{table}

\begin{figure}[ht!]\centering
\begin{tabular}{@{}ccc@{}}
\includegraphics[width=0.32\textwidth]{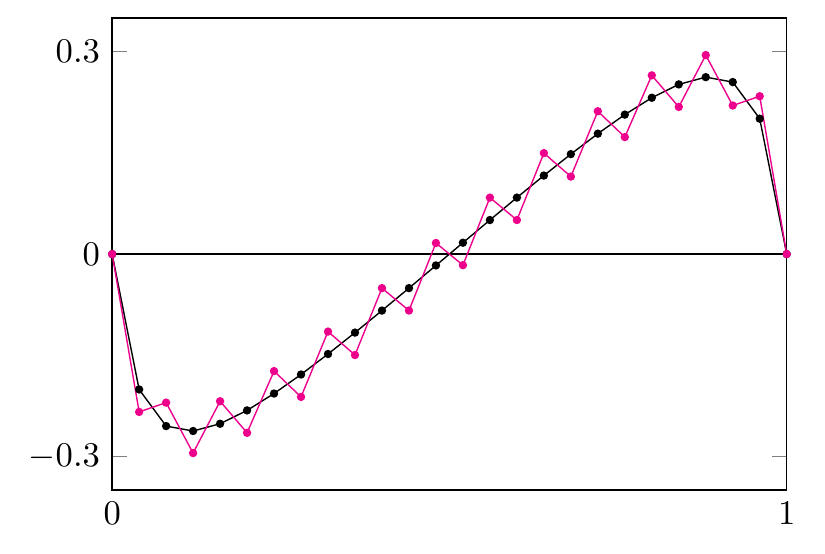} &
\includegraphics[width=0.32\textwidth]{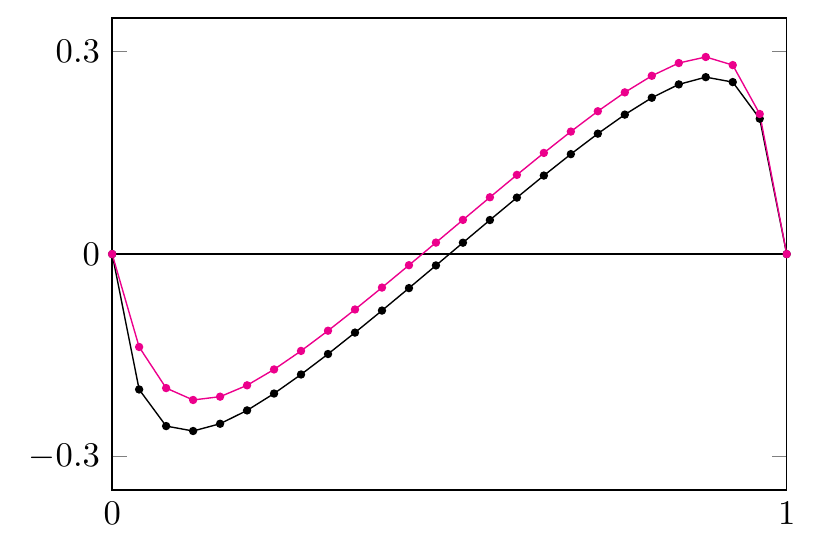} &
\includegraphics[width=0.32\textwidth]{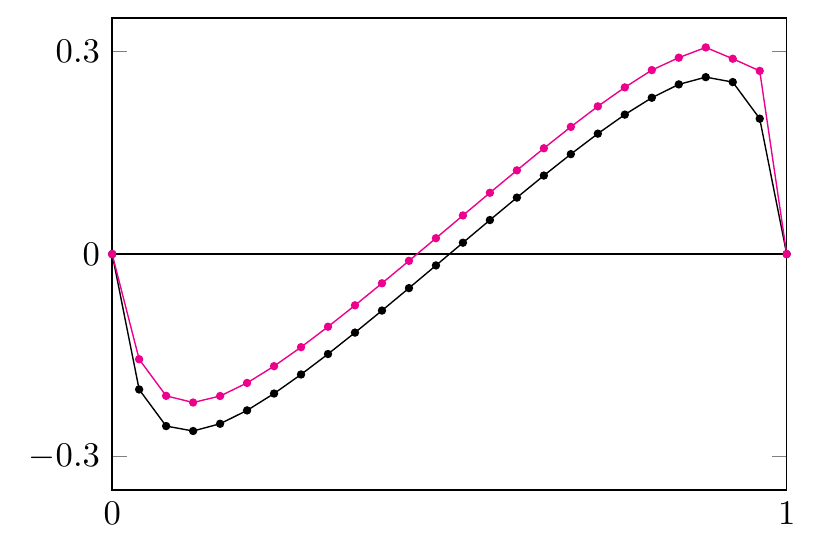}\\
\end{tabular}
\caption[]{Exact (\exactsol) and approximate (\numericalsol)
solutions for benchmark~1: centered (left), strong (middle), and
weak (right).} \label{fig:steady_solution_rough}
\end{figure}

To highlight the interest of considering upwind schemes in order to
improve the stability for large P\'eclet number, we present the
extreme situation of a pure steady-state convective problem
$\mathfrak E[\phi]=0$ with $u=1$ and $\kappa=0$ ($Pe=+\infty$) ---
benchmark~1. The manufactured solution is given by
\begin{equation}
\label{eq:delta_solution}
\phi(x;\delta)=\frac{1}{\pi}\left(1-\frac{2}{\pi}\arccos((1-\delta)\sin(\frac{\pi}{2}(2x-1)))\right)\left(\arctan(\frac{1}{\delta}\sin(\pi
x))\right),
\end{equation}
where parameter $\delta$ controls the roughness of the function. In
the present case, we take $\delta=0.1$ and carry out the numerical
simulation with $I=25$. We display in
Figure~\ref{fig:steady_solution_rough} the shape of the solution and
report oscillations for the centered scheme approximation while the
weak and strong upwind schemes eliminate the numerical artefact.

\section{Time-dependent convection-diffusion equation}

We now turn to the time-dependent problem, considering the
one-dimensional, $1$-periodic in space, convection-diffusion
equation. We seek function $\phi=\phi(x,t)$ solution of
\begin{equation}
\phi_t=\mathfrak E[\phi]+f, \quad\text{in } \Omega\times
(0,t_{\text{f}}], \label{UnstationaryConvDiff1DEquation}
\end{equation}
where $f=f(x,t)$ is a regular, $1$-periodic in space, source term,
and $t_{\text{f}}>0$ is the final time. Initial condition is
prescribed with $\phi(x,0)=\phi^0(x)$, $x\in\Omega$, while the
periodic condition reads $\phi(0,t)=\phi(1,t)$, $t>0$.

Applying the method of lines, we seek for an approximation of the
solution of the ordinary differential system of equations
\begin{equation}
\label{eq:difsystem} \frac{\text{d}\Phi}{\text{d}t}=A\Phi+F,
\quad\text{in } (0,t_{\text{f}}],
\end{equation}
where vector $\Phi=\Phi(t)=(\phi_i(t))_{i=1}^I$,
$\phi_i(t)\approx\phi(x_i,t)$, and $A=A(\theta_3,\theta_4,Pe)$ is
the circulant matrix associated to the scheme parameterised with
$\theta_3$, $\theta_4$, and $Pe$. On the other hand,
$F=F(t)=(f_i(t))_{i=1}^I$, $f_i(t)=f(x_i,t)$, while the initial
condition is given by $\phi_i(0)=\phi(x_i,0)$.

\subsection{Time discretisation}

We aim at designing a multi-stage Runge-Kutta (RK) method to compute
numerical approximations in time that provides the better trade-off
between stability, accuracy, and positivity taking into account the
cell P\'eclet number $Pe$ and the parameterisation of the spatial
scheme with respect to $\theta_3$ and $\theta_4$.

Let $N$ be a positive integer and $(t^n)_{n=0}^N$ the discrete
times. We consider the uniform subdivision $t^n=n\Delta t$,
$n=0,\ldots,N$, with the time step $\Delta
t=\frac{t_{\text{f}}}{N}$. The generic $s$-stage Runge-Kutta method
to solve the initial value ODE system~\eqref{eq:difsystem} is given
by
\begin{align*}
&\Phi^{n,j}=\Phi^n+\Delta t\sum_{\ell=1}^{s}a_{j\ell}\mathcal K^{n,\ell},\\
&\Phi^{n+1}=\Phi^{n}+\Delta t \sum_{j=1}^{s}b_j \mathcal K^{n,j},\\
&\mathcal K^{n,j}=A \Phi^{n,j}+F(t^{n,j}), \quad j=1,\ldots,s,
\end{align*}
where $\Phi^{n}=(\phi_i^{n})_{i=1}^I$,  $\phi_i^{n}\approx
\phi(x_i,t^{n})$, and $\Phi^{n,j}=(\phi_i^{n,j})_{i=1}^I$,
$\phi_i^{n,j}\approx \phi(x_i,t^{n,j})$, with $t^{n,j}=t^n+\Delta t
c_j$, are the intermediate time sub-steps.

We store the entries $(a_{j\ell})$, $(b_j)$, and $(c_j)$ in matrix
$A_\text{BT} \in\mathbb R^{s\times s}$, vectors
$b_\text{BT}\in\mathbb R^s$ and $c_\text{BT}\in \mathbb R^s$
respectively, presented in a table called \textit{Butcher tableau}:
\begin{center}
\begin{minipage}{.1\textwidth}
\begin{tabular}{@{}c|c@{}}
$c_\text{BT}$ & $A_\text{BT}$ \\
\hline
& $b_\text{BT}$\\
\end{tabular}
\end{minipage}
\begin{minipage}{0.03\textwidth}
=
\end{minipage}
\begin{minipage}{.15\textwidth}
\begin{tabular}{@{}c|ccc@{}}
$c_1$ & $a_{11}$ & $\cdots$ & $a_{1s}$ \\
$\vdots$ & $\vdots$ & & $\vdots$ \\
$c_s$& $a_{s1}$ & $\cdots$ & $a_{ss}$ \\
\hline
& $b_1$ & $\cdots$ & $b_s$\\
\end{tabular}\,.
\end{minipage}
\end{center}
Notice that the explicit Runge-Kutta is achieved if $a_{ij}=0$ for
$i\leq j$.

Equation~\eqref{eq:difsystem} leads to the uncoupled linear
differential system
\begin{equation}
\label{edo_eigenvalue} \frac{\text{d}\tilde
\phi_i}{\text{d}t}=\lambda_i(\theta_3,\theta_4,Pe)\tilde
\phi_i+\tilde f_i, \qquad \tilde \phi_i(0)=\tilde \phi_i^0, \quad
i=1,\ldots,I,
\end{equation}
taking into account the circulant matrix $A$ eigenvalues given by
equation~\eqref{eigenvalues_A}, where $\tilde \phi_i=\tilde
\phi_i(t)$ and $\tilde f_i=\tilde f_i(t)$ are the projections of
$\phi_i$ and $f_i$, respectively, on the eigenbasis.

To deal with the stability of the Runge-Kutta scheme, we consider
the homogeneous problem deriving from \eqref{edo_eigenvalue}, by
cancelling the source term $\tilde f_i$. Let $z_i=\Delta t
\lambda_i(\theta_3,\theta_4,Pe)$. The $s$-stage order $p$ explicit
Runge-Kutta scheme for the homogeneous equation associated with
equation~\eqref{edo_eigenvalue} reads
\begin{equation}
\label{stabpoly} {\tilde
\phi}^{n+1}_i=R_{ps}(z_i;w_{p+1},\ldots,w_{s}){\tilde \phi}^{n}_i,
\quad \textrm{with} \quad
R_{ps}(z;w_{p+1},\ldots,w_{s})=1+\sum_{k=1}^{p}\frac{z^k}{k!}+\sum_{k=p+1}^{s}w_kz^k,\quad
z\in\mathbb C,
\end{equation}
where $R_{ps}$ is the polynomial transfer function and
$w_k\in\mathbb R$, $k=p+1,\ldots,s$, stand for the free parameters
when $p<s$. No free parameters are available if $s=p$ and we just
denote $R_p(z)\equiv R_{pp}(z;)$ for the sake of simplicity. We
recall that stability scheme is achieved when the complex values $z$
are such that the absolute value of polynomial $R_{ps}$ is lower
than one and the set
\begin{equation}
\label{stabregion} \SR_{ps}(w_{p+1},\ldots,w_{s})=\{z \in
\mathbb{C}: |R_{ps}(z;w_{p+1},\ldots,w_{s})|< 1\},
\end{equation}
characterises the absolute stability region~\cite{Hunds}. In
addition, for $s=p$, we denote $\SR_p\equiv \SR_{pp}$.

\begin{figure}[ht]
\centering
\includegraphics{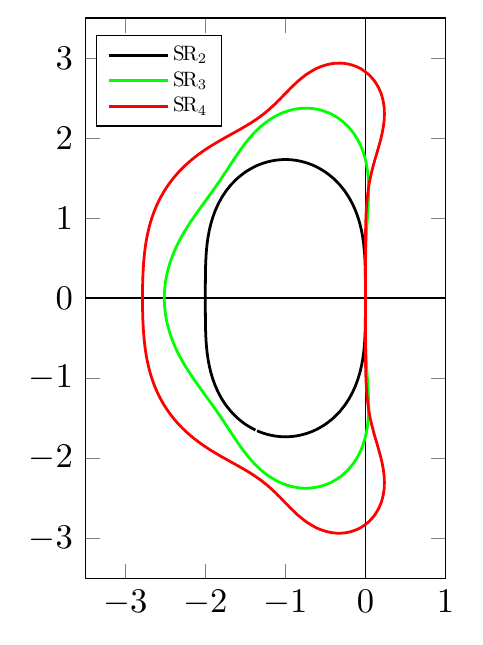}
\caption{Stability regions of classical Runge-Kutta methods
RK$_{2}$, RK$_{3}$, and RK$_{4}$.} \label{stability234}
\end{figure}

As an example, we plot in Figure~\ref{stability234} the stability
regions for three classical schemes: (i) $\SR_{2}$ --- order 2 with
$w_3=w_4=0$, (ii) $\SR_{3}$ --- order 3 with $w_3=\frac{1}{6}$,
$w_4=0$, and (iii) $\SR_{4}$ --- order 4 with $w_3=\frac{1}{6}$,
$w_4=\frac{1}{24}$.

\subsection{Design of Optimal Absolute stability regions}
The goal is to design an efficient explicit Runge-Kutta method
allowing step sizes as large as possible, still preserving the
linear stability. The key idea of the optimisation of Runge-Kutta
methods for temporal integration is to determine the parameters that
yield both the largest stability limit and the highest accuracy
(expressed in terms of dissipation and dispersion). In \cite{Ketch},
the authors focus on constructing a stability polynomial which
allows the largest absolutely stable step size and corresponding
Runge - Kutta method (number of stages and Butcher tableau) for a
given problem when the spectrum of the initial value problem is
known. They formulated a stability optimisation problem and
constructed an algorithm based on convex optimisation techniques.
They considered a global convergence in the case that the order of
approximation is one and the spectrum encloses a star-like region.
Their optimality criterium is the stability of the method.
Ait-Haddou in \cite{Ait-Haddou} used the theory of polar forms of
polynomials to obtain sharp bounds on the radius of the largest disc
(absolute stability radius), and on the length of the largest
possible real interval (parabolic stability radius) to be inscribed
in the stability region associated to the stability polynomial of an
explicit Runge-Kutta method.  Schlegel et. al. \cite{Schlegel}
constructed a multirate time-step integration method for the
convection equation. The method decouples different physical regions
so that the time step size constraint becomes a local instead of a
global restriction. Moreover, Schlegel introduced a generic
recursive multirate Runge-Kutta scheme of third order accuracy. In
\cite{Krivovichev}, the author developed an optimization of the
explicit two-derivative sixth-order Runge–Kutta method in order to
obtain low dissipation and dispersion errors. The method depends on
two free parameters, used for the optimisation and the spatial
derivatives are discretized by finite differences and
Petrov–Galerkin approximations. \vskip 1em

In this study, we focus on 4-stage RK schemes ($s=4$) at least
second-order ($p\geq 2$) as a guideline for a more general
situation. As indicated by Table~\ref{fig:centered_weak_upwind}, the
shape of spectral curve $\rks\to\lambda(\rks;\theta_3,\theta_4,Pe)$
highly depends on the cell P\'eclet number and the parameters
values. On the other hand, stability regions are controlled by $w_3$
and $w_4$ parameters and should be adapted in function of the
spectrum curve to optimally embedded the curve into
$\SR_{ps}(w_{3},w_{4})$. Such an optimisation problem is almost
intractable due to the high non-linearity involved in the
construction of the functional to minimise and we observe there
exist three major scenarios: (A) the spectral curve is almost
vertical, (B) almost horizontal, and (C) an intermediate case (see
Table~\ref{fig:centered_weak_upwind}). Consequently, we aim to
determine two absolute stability regions corresponding to the two
extreme scenarios.

\subsubsection{The imaginary axis}
Scenario (A) takes place for low diffusion schemes where the
spectral curve is getting closer to the vertical axis as long as the
P\'eclet number increases. We also deal with the fourth-order
centered scheme where the spectrum lies ``near'' the left of the
imaginary axis. Consequently, one has to design a RK scheme by
seeking real constants $w_3$ and $w_4$ such that the stability
region includes the largest segment of the imaginary axis centred at
the origin.

The solution is given by the following optimisation problem
\[
\max_{w_3,w_4\in\mathbb R} \, \Big \{\eta;\
[-\eta\mathfrak{i},\eta\mathfrak{i}]\subset \SR_{24}(w_3,w_4)\Big
\},
\]
but no analytical solution can be exhibited. Nevertheless, in
\cite{Hunds}, the authors present a solution of the optimal problem
when one maximises $\eta$ for  the three-parameters functional
$\SR_{14}(w_2,w_3,w_4)$ (which contains the particular case
$\SR_{24}(w_3,w_4)$, since we have one more free parameter). It is
shown that polynomial
\[
P_4(z)=1+z+\frac{5}{9}z^2+\frac{4}{27}z^3+\frac{4}{81}z^4
\]
provides the best largest segment of the imaginary axis with
$\eta_{\max}=3$.

\begin{figure}[ht]
\centering
\includegraphics{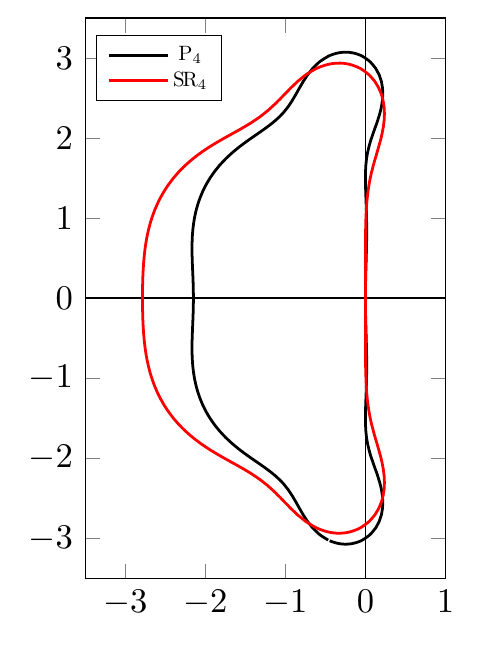}
\caption{Stability region associated to polynomials $P_4$ and
R$_{4}$.} \label{fig:P4_RK44}
\end{figure}

We plot in Figure~\ref{fig:P4_RK44} the two stability regions
associated to the popular RK$_4$ scheme and the optimal solution
proposed in~\cite{Hunds}. We observe that the RK$_4$ scheme provides
an excellent approximation with $\eta_{\max}=2\sqrt{2}$,
\cite{Hunds}. Of course, the optimal case would provide a slightly
bigger $\eta_{\max}$, but we consider that, for our application, the
RK$_4$ is an excellent candidate for the first scenario.

\subsubsection{The real axis}

\begin{figure}[ht]
\centering
\includegraphics{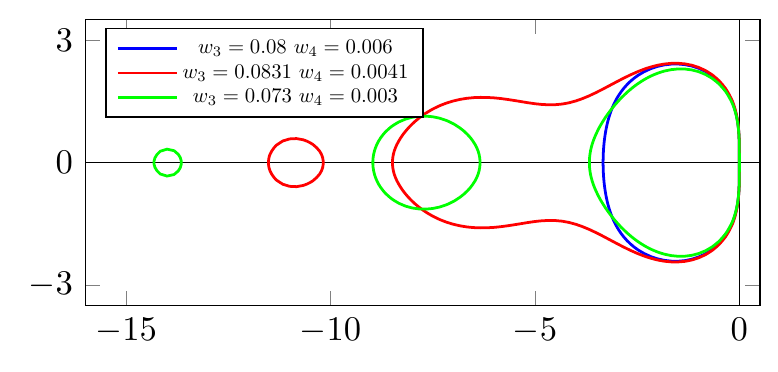}
\caption{Stability regions for different values of $w_3$ and $w_4$.}
\label{fig:1-stability}
\end{figure}

Scenario (B) concerns numerical schemes with large diffusion
characterised by low P\'eclet numbers. Therefore, we seek real
constants $w_3$ and $w_4$ such that the stability region includes
the largest segment of the negative real axis starting at the
origin. An additional difficulty is that the stability region may
not be a simply connected region as exemplified in
Figure~\ref{fig:1-stability}. Consequently, we only consider the
first connected region which contains the origin as the effective
stability region. There are also stabilized explicit Runge - Kutta
methods as, for example, Runge - Kutta - Chebyshev methods (RKC)
dedicated to extended real stability intervals and useful for
semi-discrete parabolic problems. A second-order RKC method was
initially proposed by van der Houwen and Sommeijer \cite{Houwen} and
a family of second- and fourth- order Orthogonal - Runge - Kutta -
Chebyshev methods (ROCK) were proposed by Abdulle and Abdulle and
Medovikov, \cite{Abdulle 2, Abdulle Medo} but as we said these
methods are based on Runge - Kutta and for our purposes, we only
deal with original Runge - kutta.

Riha, in \cite{Riha}, proved that the optimal stability polynomial
$\bar R_{ps}(z)$ of order $p$ and degree $s$
\[
\bar R_{ps}(z; \bar w_1,\ldots,\bar w_s)=1+\sum_{k=1}^p
\frac{z^k}{k!}+\sum_{k=p+1}^s \bar w_kz^k
\]
such that
\[
|\bar R_{ps}(z; \bar w_1,\ldots,\bar w_s)|\leq 1, \quad z\in
[-\zeta,0],
\]
is unique and satisfies the so-called ripple property:
\begin{prop}[ripple property]
Polynomial $R_{ps}(z;w_1,\ldots,w_s)$ satisfies the ripple property,
if and only if, there exist $s-p+1$ points $x_0<x_1<\ldots
<x_{s-p}<0$, with $x_0=-\zeta_{\max}$, such that
\begin{alignat*}{3}
& R_{ps}(x_i;w_1,\ldots,w_s) &&= -R_{ps}(x_{i+1};w_1,\ldots,w_s),  && \quad i=0,1,\ldots, s-p-1,\\
& |R_{ps}(x_i;w_1,\ldots,w_s)| &&= 1, && \quad i=0,1,\ldots, s-p.
\end{alignat*}
\end{prop}
There are no explicit analytic expressions for the optimal $s-p+1$
coefficients $w_k$, $k=p+1,\ldots,s$, but the ripple property has
been used to construct approximations to the optimal stability
polynomial $R_{ps}(z)$ as proposed for example in~\cite{Hunds,
Skorv}, where the authors show that the optimal bound $\zeta_{\max}$
depends on $s$ and satisfies $\zeta_{\max}=c_ps^2$, $c_p \in
\mathbb{R}$, asymptotically with $s \to \infty$.

\begin{figure}[ht]
\centering
\includegraphics[width=.45\linewidth]{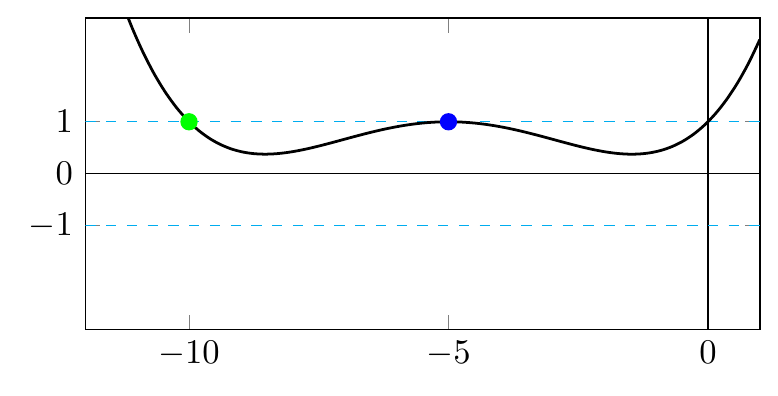}
  \includegraphics[width=.45\linewidth]{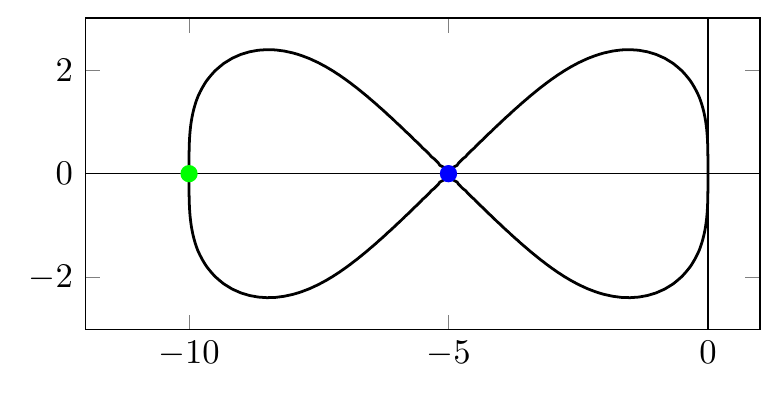}
\caption{Baker polynomial $B_4=R_{24}(x;\frac{2}{25},
\frac{1}{250})$ graph (left) and the corresponding stability region
(right).} \label{fig:Baker_polynomial}
\end{figure}

For $p=2$, there is a suitable approximate polynomial $B_s$ based on
Chebyshev polynomials, given by M. Bakker in 1971, that generates
about $80\%$ of the optimal interval and for $s=4$, the Bakker
polynomial reads
\begin{equation}
\label{Bakker}
B_4(z)=1+z+\frac{z^2}{2}+\frac{2z^3}{25}+\frac{z^4}{250},
\end{equation}
where $\zeta_{\max}=10$, $w_3=\frac{2}{25}$, and
$w_4=\frac{1}{250}$. We plot in Figure~\ref{fig:Baker_polynomial}
the Baker polynomial representation and the corresponding stability
region.

\begin{figure}[ht]
\centering
  \includegraphics[width=.45\linewidth]{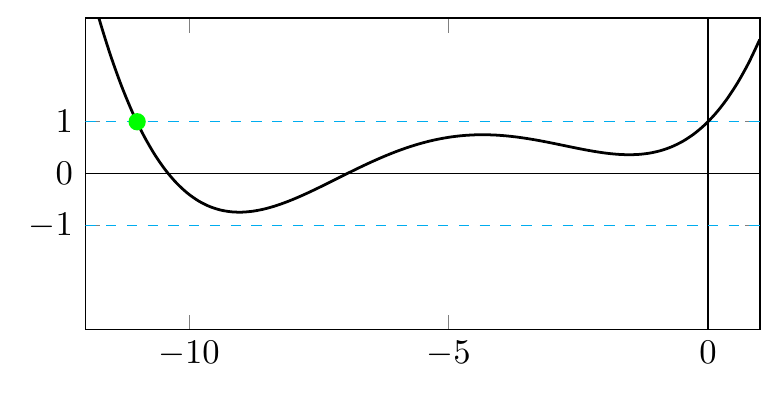}
  \includegraphics[width=.45\linewidth]{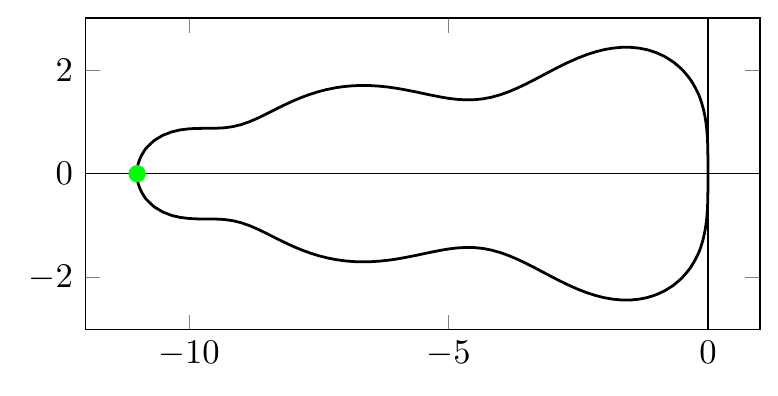}
\caption{ Polynomial $R_{24}\left(x;0.083, 0.0042\right)$ graph
(left) and the corresponding stability region (right). }
\label{fig:low_pe_polynomial}
\end{figure}

Such a stability region is suitable for $Pe=0$ but not for small
values of $Pe$ since it is not radial at the origin. Consequently,
we have considered a small perturbation of parameters $w_3$ and
$w_4$ in order to produce a radial region still preserving a large
interval on the real negative axis. We found a good trade-off with
$w_3=0.0834$ and $w_4=0.0042$ providing $\zeta_{\max}\approx11$. We
plot in Figure~\ref{fig:low_pe_polynomial} the polynomial curve and
the corresponding region.

\subsubsection{Design of a non-negative scheme}

A critical issue in some problems is to guarantee the positivity of
the solution. Indeed, temperature, concentration, or density are
physical quantities that must be non negative at the discrete level.
Figure~\ref{fig:low_pe_polynomial} shows that the polynomial does
not satisfy this criterion. For instance, if $x \in
]-10.3997,-6.9570[$, the scheme is stable but the polynomial is
negative leading to a sign change.

We then consider a more constrained problem and define the new
optimisation problem
\begin{equation*}
\displaystyle\max_{w_3,w_4\in\mathbb R} \Big \{ \zeta; 0.01\leq
R_{24}(x;w_3,w_4)\leq 0.7, \forall x\in[-\zeta,0]\Big \},
\end{equation*}
to create the appropriated region. The lower bound $0.01$ and upper
bound $0.7$ are prescribed in order to provide a radial domain.
Indeed, extreme bound values $0.0$ and $1.0$ provide a stability
region similar to Bakker polynomial displayed in
Figure~\ref{fig:Baker_polynomial}.

\begin{figure}[ht]
\centering
\includegraphics[width=.45\linewidth]{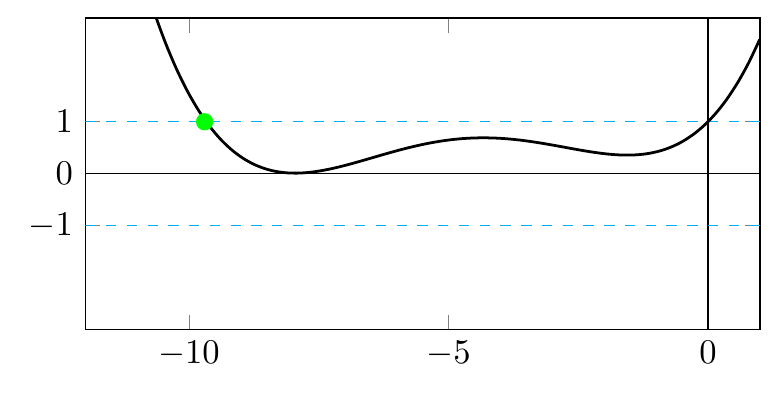}
\includegraphics[width=.45\linewidth]{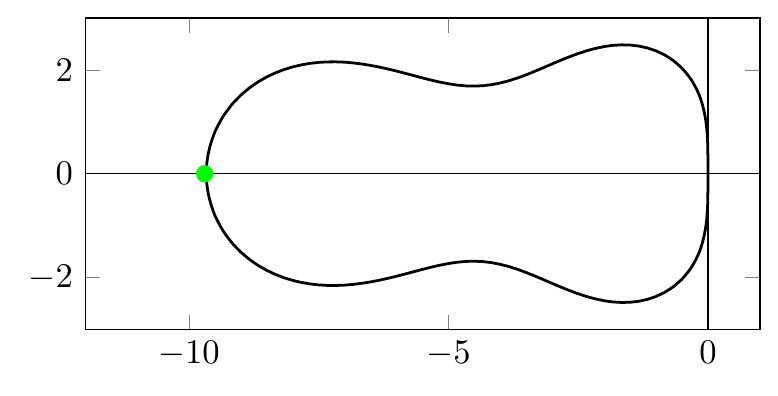}
\caption{ Polynomial $R_\text{D}(x)$ graph (left) and the
corresponding stability region (right).}
\label{fig:low_pe_positive_polynomial}
\end{figure}

Numerical solution of the optimal problem provides
$w_3=\frac{603}{6998}$ and $w_4=\frac{15}{3212}$ with
$\zeta_{\max}\approx9.43$ and we plot in
Figure~\ref{fig:low_pe_positive_polynomial} the polynomial curve and
the corresponding stability region. We obtain a large radial
stability region with $\zeta_{\max}$ close to the one provided by
the Bakker polynomial. On the other hand, the domain still preserves
an important part of the horizontal axis in comparison to the non
positive optimal case given in Figure~\ref{fig:low_pe_polynomial}.

To conclude the analysis, scenario (A) concerns the low diffusive
situation and the $R_4$ polynomial is well adapted for such
situation. On the other hand,  scenario (B) deals with high
diffusive operator and the polynomial
\[
R_\text{D}(z)=R_{24}\left(z;\frac{603}{6998}, \frac{15}{3212}\right)
\]
provides an excellent stability region.

\subsection{Butcher tableau construction}

Stability issue provides coefficients $w_3$ and $w_4$ for polynomial
$R_\text{D}$ and one has to compute corresponding entries of Butcher
tableau. We briefly outline the method for the general case
$R_{24}(z;w_3,w_4)$ and provide the Butcher's tableau for the two
scenarios, taking into account a crucial restriction: the two
tableaux must have the same time sub-steps for the sake of
compatibility.

An explicit 4-stage Runge-Kutta method is given by the generic
Butcher tableau
\begin{center}
\begin{tabular}{@{}c|ccccc@{}}
$c_1$ & 0 & 0 & 0 & 0 \\
$c_2$ & $a_{21}$ & 0 & 0 & 0 \\
$c_3$ & $a_{31}$ & $a_{32}$ & 0 & 0 \\
$c_4$ & $a_{41}$ & $a_{42}$ & $a_{43}$ & 0 \\\hline & $b_1$ & $b_2$
& $b_3$ & $b_4$
\end{tabular}
\end{center}

Consistency yields
\[
c_i=\sum_{j=1}^{i-1}a_{ij}, \quad i=1, \ldots, 4,
\]
while the second-order assumption provides the constraints
\begin{equation}
\label{cond_2nd_order}
\begin{bmatrix} 1 & 1 & 1 & 1 \end{bmatrix} b=1, \quad b^\text{T}c=\frac{1}{2}.
\end{equation}

Additional constraints are considered:
\begin{itemize}
\item The Runge-Kutta scheme is fourth-order for the non-homogeneous equation $\phi'(t)=f(t)$,
\begin{equation}
\label{cond_non_hom}
\begin{bmatrix}(c_1)^2 & (c_2)^2 & (c_3)^2 & (c_4)^2\end{bmatrix} b=\frac{1}{3}, \qquad \begin{bmatrix}(c_1)^3 & (c_2)^3 & (c_3)^3 & (c_4)^3\end{bmatrix} b=\frac{1}{4}.
\end{equation}
\item The four stage method has to satisfy the consistency condition to suit the polynomial $R_{24}(z;w_3,w_4)$,
\begin{equation}
\label{cond_hom} b^\text{T}Ac=w_3, \qquad b^\text{T}A^2c=w_4.
\end{equation}
\end{itemize}

Conditions~\eqref{cond_2nd_order}, \eqref{cond_non_hom}, and
\eqref{cond_hom} are equivalent to the system
\begin{equation}
\label{vandermonde} V^\text{T}b^\text{T}=
\begin{bmatrix}
1 & 1 & 1 & 1 \\
c_1 & c_2 & c_3 & c_4 \\
c_1^2 & c_2^2 & c_3^2 & c_4^2 \\
c_1^3 & c_2^3 & c_3^3 & c_4^3 \\
\end{bmatrix}
\begin{bmatrix}
b_1 \\
b_2 \\
b_3 \\
b_4 \\
\end{bmatrix}
=\begin{bmatrix}
1 \\
1/2 \\
1/3 \\
1/4 \\
\end{bmatrix}
\end{equation}
together with equations
\begin{align}
& b_3a_{32}c_2+b_4(a_{42}c_2+a_{43}c_3)=w_3 \label{system a_b_1}\\
& b_4a_{43}a_{32}c_2=w_4 \label{system a_b_2}.
\end{align}
Assuming $c_1=0$, the determinant of the Vandermonde matrix $V$ is
$(c_4-c_3)(c_3-c_2)c_2$ and two situations arise:
\begin{itemize}
\item If $c_2 \neq 0$, $c_4 \neq c_3$, and $c_3 \neq c_2$, the system has a unique solution $b$. Hence in this case, we choose the free parameters to be $c_2$, $c_3$, $c_4$, and $a_{43}$. From the Vandermonde system, we obtain $b$.
\item If $c_2 = 0$ or $c_4 = c_3$ or $c_3 = c_2$, we have a dependent linear system. For instance, assuming that $c_2=c_3$, system~\eqref{vandermonde} turns to be
\begin{equation}
\label{vandermonde1}
\begin{bmatrix}
1 & 1 & 1 & 1 \\
0 & c_2 & c_2 & c_4 \\
0 & 0 & 0 & c_4(c_4-c_2) \\
0 & 0 & 0 & 0 \\
\end{bmatrix}
\begin{bmatrix}
b_1 \\
b_2 \\
b_3 \\
b_4 \\
\end{bmatrix}
=
\begin{bmatrix}
1 \\
1/2 \\
1/3-c_2/2 \\
1/4-(c_4+c_2)/3+c_4c_2/2 \\
\end{bmatrix}.
\end{equation}
The system has several degrees of freedom we have to fix with
relations
\[
\frac{1}{4}-\frac{c_4+c_2}{3}+\frac{c_4c_2}{2}=0, \quad
c_4(c_4-c_2)\neq 0, \quad c_2\neq 0.
\]
We get vector $b$ with
\[
b_4=\frac{2-3c_2}{6c_4(c_4-c_2)},\quad
b_3=\frac{3c_4-2}{6c_2(c_4-c_2)}-b_2,\quad
b_1=\frac{c_4+c_2-1}{6c_2c_4}.
\]
Notice that the user has to choose the free parameters $c_2$,
$a_{43}$, and $b_2$.
\end{itemize}

From vectors $c$ and $b$, we compute the remaining Butcher tableau
elements $a_{21}$, $a_{31}$, $a_{32}$, $a_{41}$, and $a_{42}$ and we
obtain
\[
a_{21} = c_2,\,\,\, a_{32} = \frac{w_4}{b_4 a_{43} c_2},\,\,\,
a_{42} = \frac{w_3-b_4 a_{43} c_3-b_3 c_2 a_{32} }{b_4c_2},\,\,\,
a_{31} = c_3-a_{32},\,\,\, a_{41} = c_4-a_{42}-a_{43}.
\]
At last, we present in Table~\ref{tab::RK4_RK_D} the two Butcher
tableaux with 4-stage corresponding to $R_{4}$ and $R_\text{D}$
polynomials. We tag the associated methods as RK$_4$ and
RK$_\text{D}$, respectively. It is important to notice that the
RK$_\text{D}$ has been designed with the same time sub-steps
according to the classical RK$_4$ scheme.

\begin{table}[ht]
\caption{Butcher tableaux for the classical scheme $\text{RK}_{4}$
(left) and RK$_\text{D}$ (right).} \label{tab::RK4_RK_D}
\begin{center}
{\notapolice
\begin{tabular}{@{}c|ccccc@{}}
0 & 0 & 0 & 0 & 0 \\
1/2 & 1/2 & 0 & 0 & 0 \\
1/2 & 0 & 1/2 & 0 & 0 \\
1 & 0 & 0 & 1 & 0 \\\hline
 & 1/6 & 1/3 & 1/3 & 1/6\\
\end{tabular}\hskip 3em
\begin{tabular}{@{}c|ccccc@{}}
0 & 0 & 0 & 0 & 0 \\
1/2 & 1/2 & 0 & 0 & 0 \\
1/2 & 334/861 & 373/3328 & 0 & 0 \\
1 & 481/3310 & 587/1655 & 1/2 & 0 \\
\hline
 & 1/6 & 0.4 & 4/15 & 1/6\\
\end{tabular}
}
\end{center}
\end{table}

\section{Optimal time step for stability}

We now reach the key point of the paper: the confrontation between
the RK stability region with the spatial operator spectrum. More
precisely, we have, on the one hand, a complex value parametric
curve  given by relation~\eqref{equ_lambda}, $\rks\in [0,1]\to
\lambda(\rks;\theta_3,\theta_4,Pe)$ that contains all the
eigenvalues $\lambda_i$ of the discrete operator. On the other hand,
applying the RK scheme with time step $\Delta t$ and setting
$z_i=\Delta t \lambda_i(\theta_3,\theta_4,Pe)$, one has to choose
$\Delta t$ small enough such that $z_i \in
\SR_{ps}(w_{p+1},\ldots,w_{s})$ to guarantee the stability.

\subsection{The CFL condition}

Let
\begin{equation}
\label{eigenCFL} z(\rks;\theta_3,\theta_4,Pe)=\Delta t
\lambda(\rks;\theta_3,\theta_4,Pe)=C_{\text{CFL}}\,\rho(\rks;\theta_3,\theta_4,Pe),
\end{equation}
where $C_{\text{CFL}}=\frac{u\Delta t}{\Delta x}$ and $\rho$ be
given by relation~\eqref{equ_rho} for $u\ne 0$. Stability condition
for the discrete problem is achieved if we satisfy the condition
\begin{equation}\label{eq:stability_CFL}
C_{\text{CFL}}\,\rho(\rks;\theta_3,\theta_4,Pe)\subset
\SR_{24}(w_{3},w_{4}),\quad \rks\in[0,1].
\end{equation}

Note that the scheme stability depend on the six parameters
$\theta_3$, $\theta_4$, $w_3$, $w_4$, $Pe$, and $C_{\text{CFL}}$ and
the goal of this section is to analyse the stability of the full
time-dependent convection-diffusion equation.
\begin{myremark}
The stability condition combines two main ingredients. On the one
hand, the shape of the stability region is characterised by function
$\rho$ that provides a complex value curve where all the eigenvalues
lie, independently of the space parameter $\Delta x$. On the other
hand, the CFL value $C_{\text{CFL}}$ scales the previous curve to
fit inside the stability region.
\end{myremark}

Assuming that scheme in space is given (parameters $\theta_3$,
$\theta_4$ are prescribed) and the scheme in time is given
(parameters $w_3$, $w_4$ are prescribed), we define the optimal CFL
curve as a function of the P\'eclet number
\[
Pe\to \widehat C_{\text{CFL}}(Pe)= \lim \sup \Big
\{C_{\text{CFL}}\geq 0, \text{ such that }
\eqref{eq:stability_CFL}\text{ holds}\Big \}.
\]
\begin{figure}[!ht]\centering
\begin{tabular}{cc}
\includegraphics{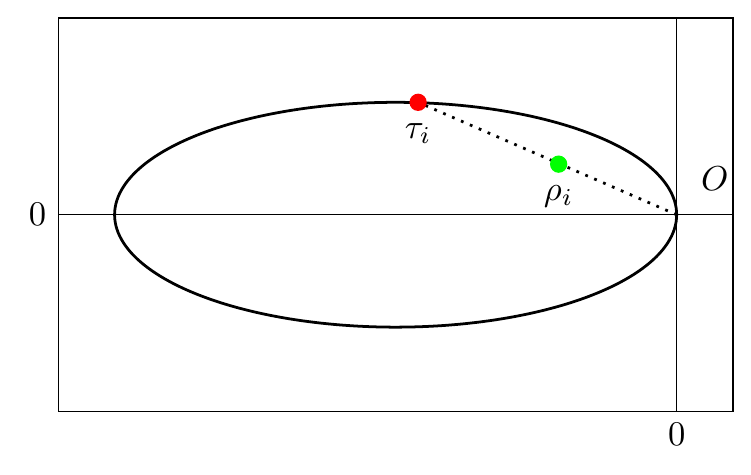} &
\includegraphics{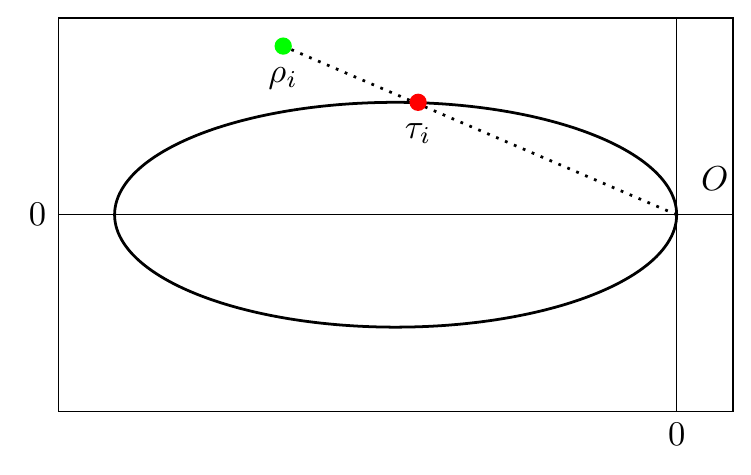}
\end{tabular}
\caption{How to calculate the limiting scaling factor: $\widehat
C_{\text{CFL}}>1$ (left) and $\widehat C_{\text{CFL}}<1$ (right).}
\label{fig::fig1}
\end{figure}
We compute $\widehat C_{\text{CFL}}(Pe)=\widehat
C_{\text{CFL}}(Pe;\theta_3,\theta_4,w_3,w_4)$ with the following
algorithm (see Figure \ref{fig::fig1}):
\begin{itemize}
\item[1.] Compute eigenvalues $\rho_i$, $i=1,\ldots,I$, from the space discretisation operator spectrum;
\item[2.] Compute intersection points $\tau_i$ of the segment line $O\rho_i$ with the stability region boundary;
\item[3.] Compute $\displaystyle\widehat C_{\text{CFL}}=\min_{i=1,\ldots,I} \frac{|\tau_i|}{|\rho_i|}$.
\end{itemize}
Afterwards the maximum stable time step is given by
\begin{equation}
    \label{delta_t_max}
\Delta t_{\max}=\frac{\widehat C_{\text{CFL}}\Delta x}{u}.
\end{equation}

\begin{myremark}
Two extreme situations require a specific treatment.
\begin{itemize}
\item If $Pe=0$, \textit{i.e.}, $u=0$, a CFL constant based on the velocity is no longer available. In that case, the spectrum is given in Table~\ref{tab:spectra_centered} for the centered scheme, Table~\ref{tab:spectra_weak} for the weak upwind, Table~\ref{tab:spectra_strong} for the strong upwind and reads
$$
\Delta t \lambda(\rks)=C_\text{CFL}\,\rho(\rks)= \frac{\kappa\Delta
t}{\Delta x^2}\, 2(\cos(2\pi \rks)-1) \left \{ \begin{array}{ll}
(\frac{7}{6}-\frac{1}{6}\cos(2\pi\rks))& \textrm{centered},\\
1&\textrm{weak},\\
(\cos(2\pi\rks)-\mathfrak{i}\sin(2\pi\rks)) &\textrm{strong}.
\end{array} \right.
$$
In this case the maximum stable time step is given by
\[
\Delta t_{\max}=\frac{\widehat C_{\text{CFL}}\Delta x^2}{\kappa}.
\]
\item If $Pe=+\infty$, \textit{i.e.}, $\kappa=0$, one has to pass to the limit to determine the spectrum curve, remaining~\eqref{delta_t_max} valid.
\end{itemize}
\end{myremark}

\begin{figure}[!ht]\centering
\begin{tabular}{@{}ccc@{}}\toprule
& centered & weak upwind\\\midrule
& $Pe=10$, $\widehat C_{\text{CFL}}=2.0935$, $\Delta t_{\max}=0.0837$ & $Pe=5$, $\widehat C_{\text{CFL}}=1.3117$, $\Delta t_{\max}=0.0525$\\
RK$_4$ & \includegraphics[valign=m]{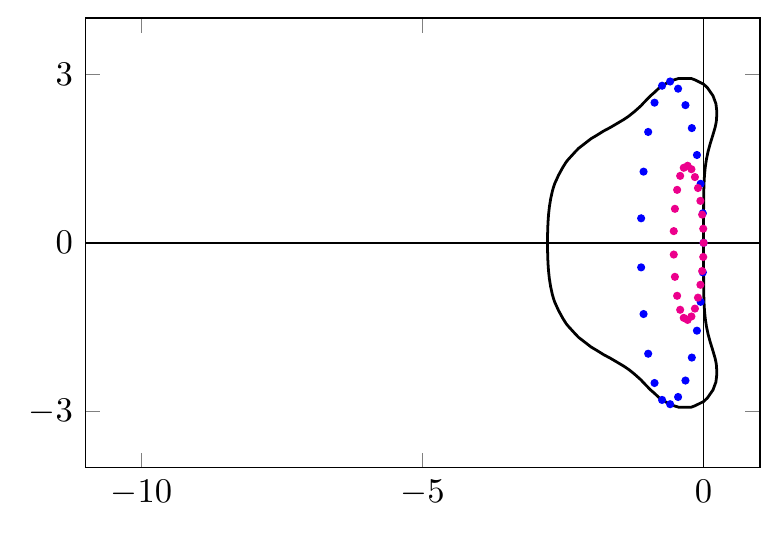}
&
\includegraphics[valign=m]{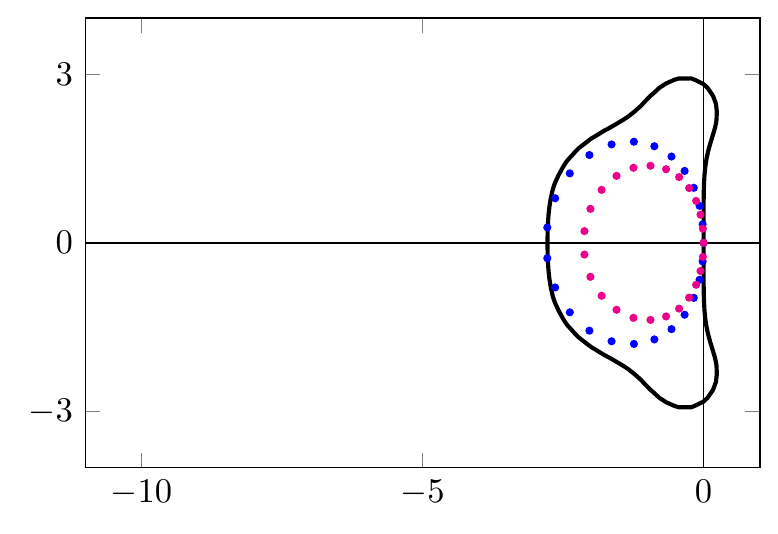} \\ \midrule
& $Pe=10$, $\widehat C_{\text{CFL}}=1.3479$, $\Delta t_{\max}=0.0539$ & $Pe=5$, $\widehat C_{\text{CFL}}=1.7948$, $\Delta t_{\max}=0.0718$\\
RK$_\text{D}$ &
\includegraphics[valign=m]{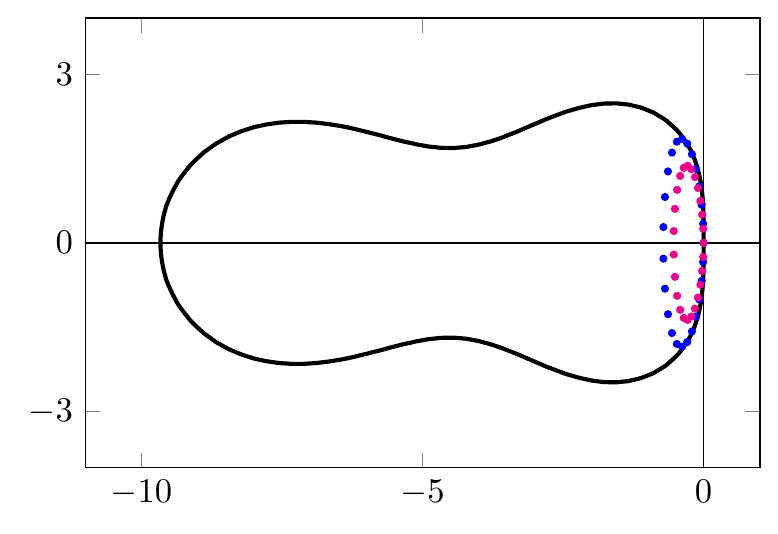} &
\includegraphics[valign=m]{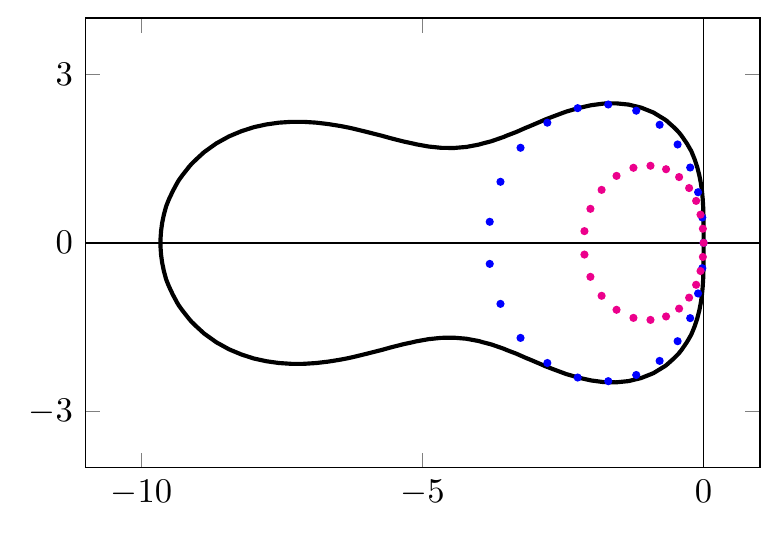}\\ \bottomrule
\end{tabular}
\caption[]{Examples of stability regions (\stabregion boundary of
stability region; \eigenwithout ``discrete'' $\rho$; \eigenwith
``discrete'' $\widehat C_{\text{CFL}}\,\rho$).}
\label{figs:stability_regions}
\end{figure}

As an example, we present in Figure~\ref{figs:stability_regions} two
situations with $Pe=5$ and $Pe=10$ where we adjust the CFL constant
to fit the spectrum into the stability region: centered scheme in
space with the RK$_4$ scheme in time (top left) and the RK$_D$ (down
left) while the pictures in second column present the weak upwind
case.

\subsubsection{Analysis of spatial schemes}

\begin{figure}[!ht]\centering
\begin{tabular}{@{}cccc@{}}\toprule
$I$ & centered & weak upwind & strong upwind\\\midrule 25 &
\includegraphics[align=c,height=4cm]{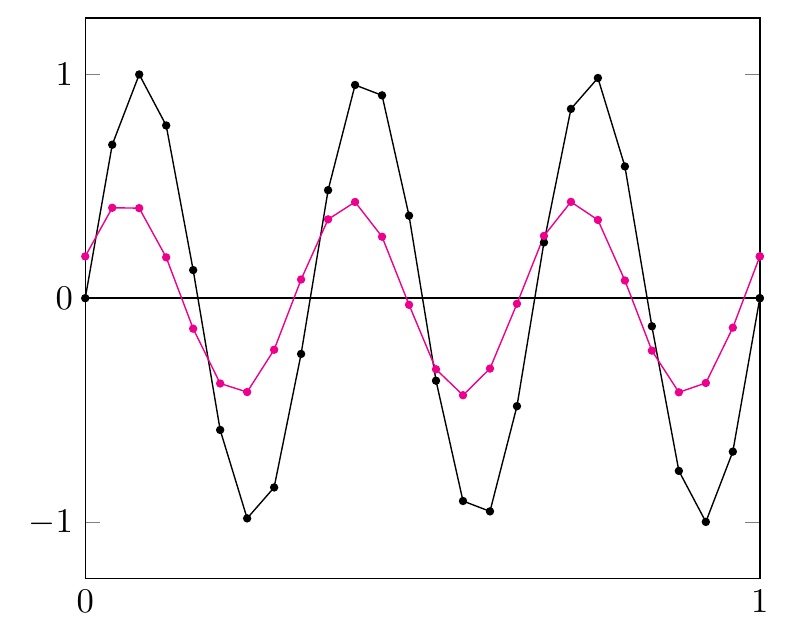} &
\includegraphics[align=c,height=4cm]{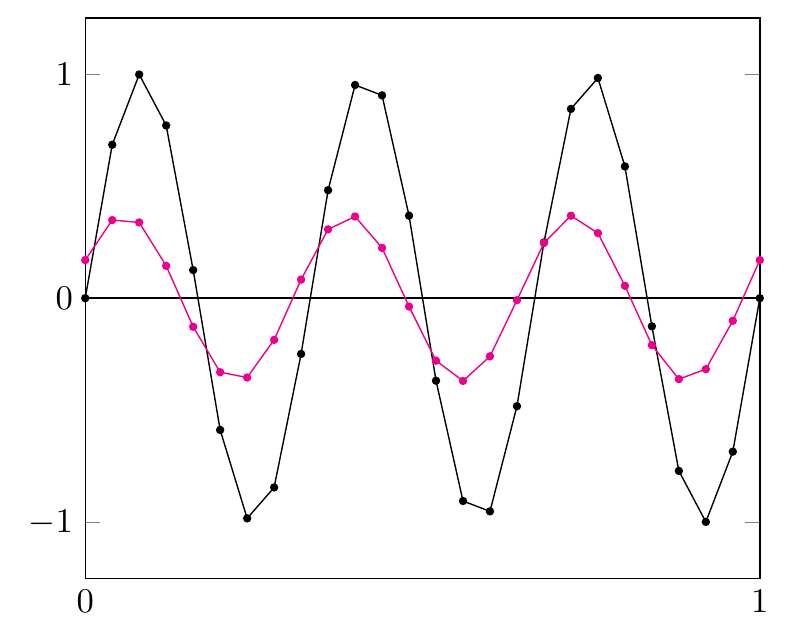} &
\includegraphics[align=c,height=4cm]{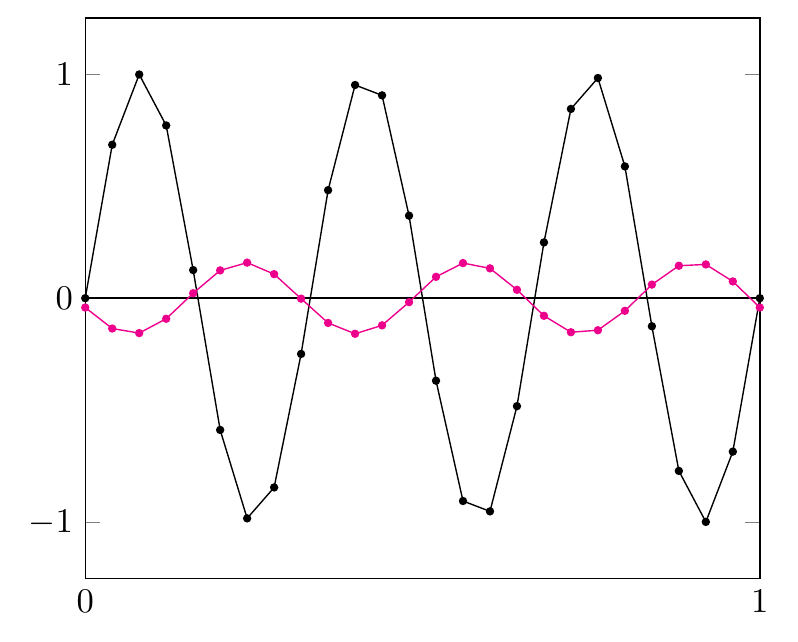} \\\midrule
50 &
\includegraphics[align=c,height=4cm]{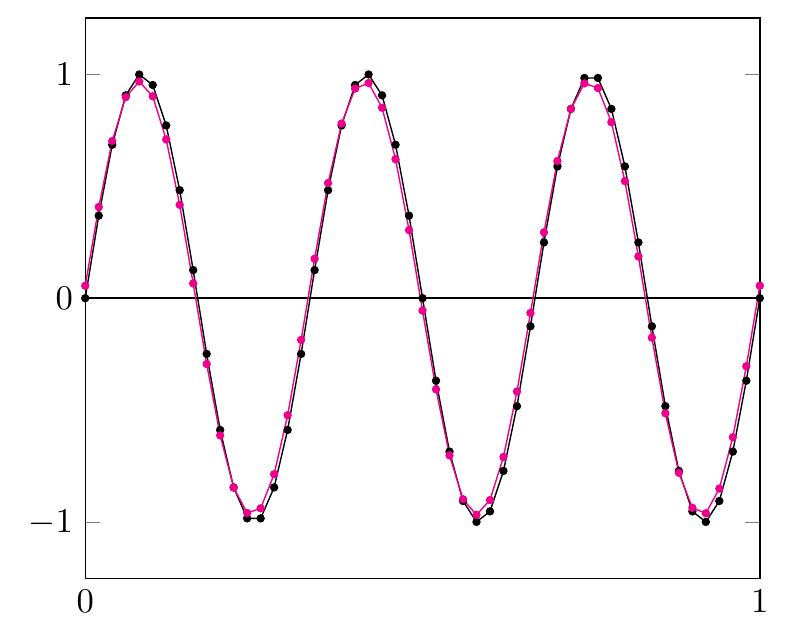} &
\includegraphics[align=c,height=4cm]{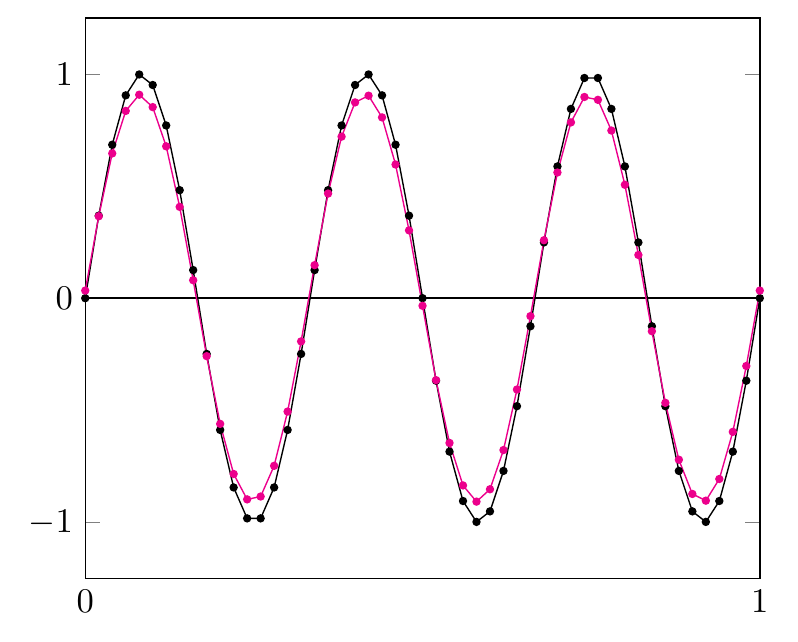} &
\includegraphics[align=c,height=4cm]{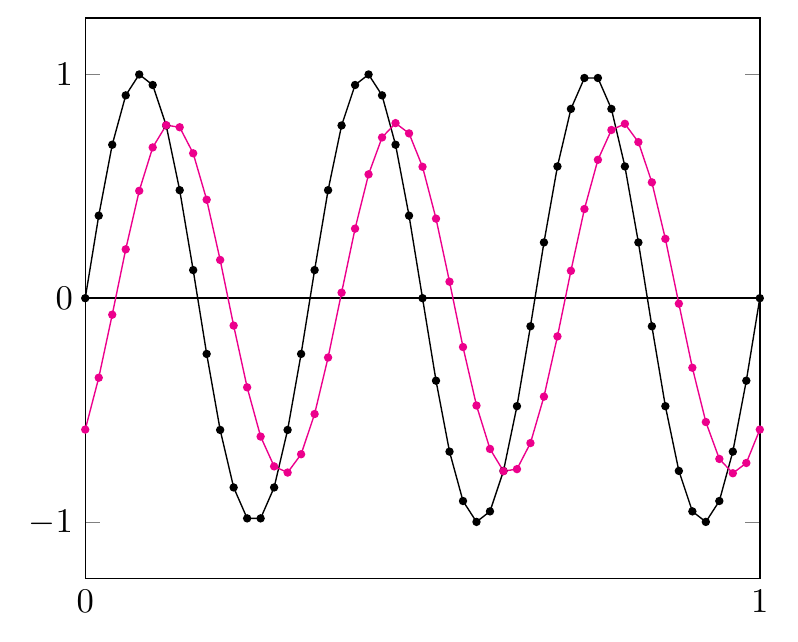} \\\bottomrule
\end{tabular}
\caption[]{Exact (\exactsol) and approximate (\numericalsol)
solutions considering RK$_4$ time scheme for benchmark 2.}
\label{fig:dispersion_strong_upwind}
\end{figure}

A full discretised scheme consists in fixing the time and spatial
scheme parameters. We have identified two schemes in time ---
RK$_\text{D}$ and RK$_4$ --- for low and high P\'eclet situations
and three schemes in space --- centered, weak upwind, and strong
upwind. We shall discard the strong upwind scheme due to the high
dispersion effect. Indeed, let us consider the advection problem
$\mathfrak E[\phi]=0$ with $u=1$ and $\kappa=0$ ($Pe=+\infty$).
Function
\[
\phi(x,t;\omega)=\sin(2\pi \omega(x-u t))
\]
with $\omega=3$, is a periodic solution. We carry out the numerical
simulation with the RK$_4$ scheme in time and the three schemes in
space with $\Delta t=\Delta t_{\max}$ (see
expression~\eqref{delta_t_max}) until the final time
$t_{\text{f}}=1$. We plot the numerical solutions in
Figure~\ref{fig:dispersion_strong_upwind} computed with $I=25$ and
$I=50$ nodes --- benchmark~2. We observe the large phase errors
produced by the strong upwind scheme which justify the choice to
discard it.

We turn to a situation where a rough function is convected ($u=1$,
$\kappa=0$, $Pe=+\infty$) and compare the centered and weak upwind
scheme. The manufactured solution is given by
\begin{equation}
\label{eq:delta_solution_xt}
\phi(x,t;\delta)=\frac{1}{\pi}\left(1-\frac{2}{\pi}\arccos((1-\delta)\sin(\pi(x-t-\frac12)))\right)\left(\arctan(\frac{1}{\delta}\sin(\pi
(x-t)))\right)
\end{equation}
with $\delta=0.01$ that corresponds to a rough solution due to the
sharp transition. The numerical solutions are computed with the
RK$_\text{4}$ scheme in time taking $\Delta t=\Delta t_{\max}$,
$I=25$, and the centered and weak upwind schemes in space ---
benchmark~3.

\begin{figure}[!ht]\centering
\begin{tabular}{@{}cc@{}}
\includegraphics[height=4cm]{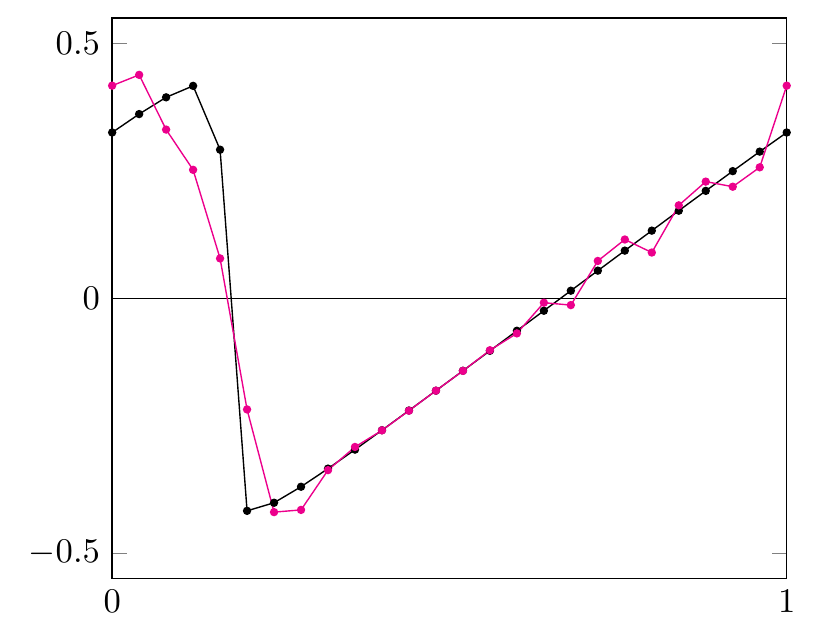} &
\includegraphics[height=4cm]{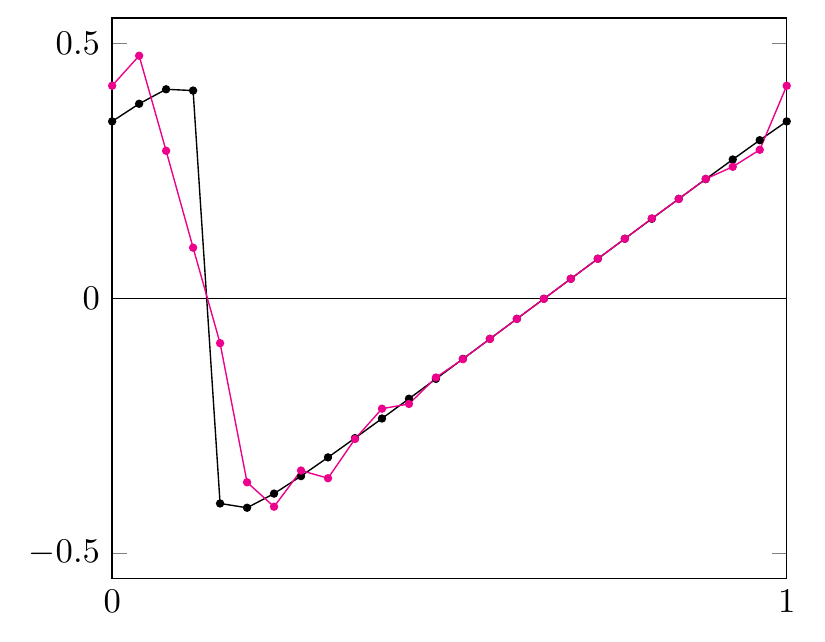}
\end{tabular}
\caption[]{Exact (\exactsol) and approximate (\numericalsol)
solutions for $t_\text{f}=2\Delta t_{\max}$ considering RK$_4$ time
scheme for benchmark 3: centered (left) and weak upwind (right).}
\label{fig:rough_solution_stability}
\end{figure}

In Figure~\ref{fig:rough_solution_stability} we present the
numerical solutions after two time steps. We observe the centered
scheme is more unstable with a larger number of over- and
undershoots. Of course, a definitive method would consist in
employing a high viscous scheme (the simple two-points upwind one)
but with a dramatic cut of the accuracy. Hence, the weak upwind
scheme is regarded as an alternative to the centered one when large
oscillations appear.

\subsubsection{Optimal CFL curves}

\begin{figure}[!ht]\centering
\begin{tabular}{@{}cc@{}}
\includegraphics[height=6cm]{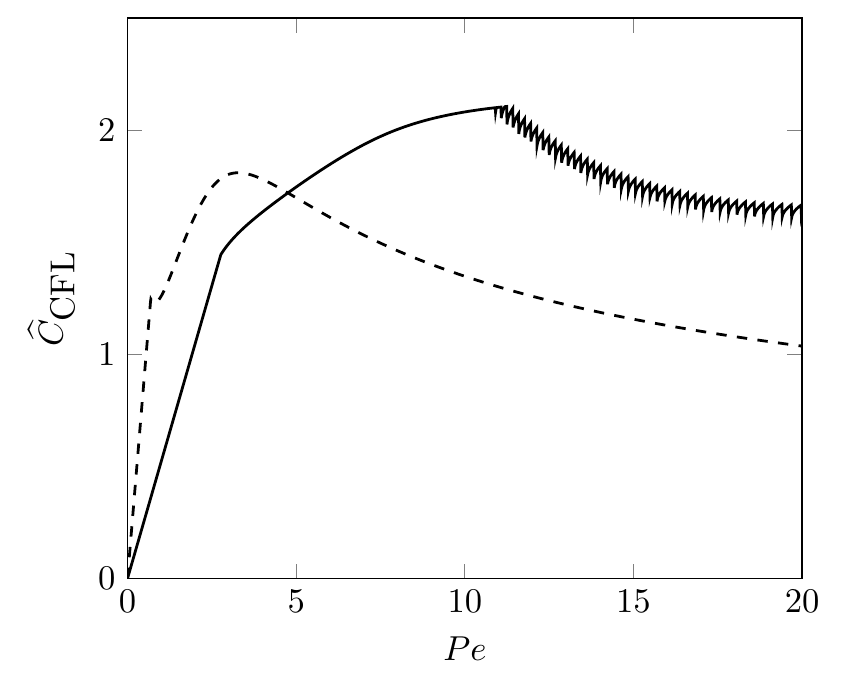} &
\includegraphics[height=6cm]{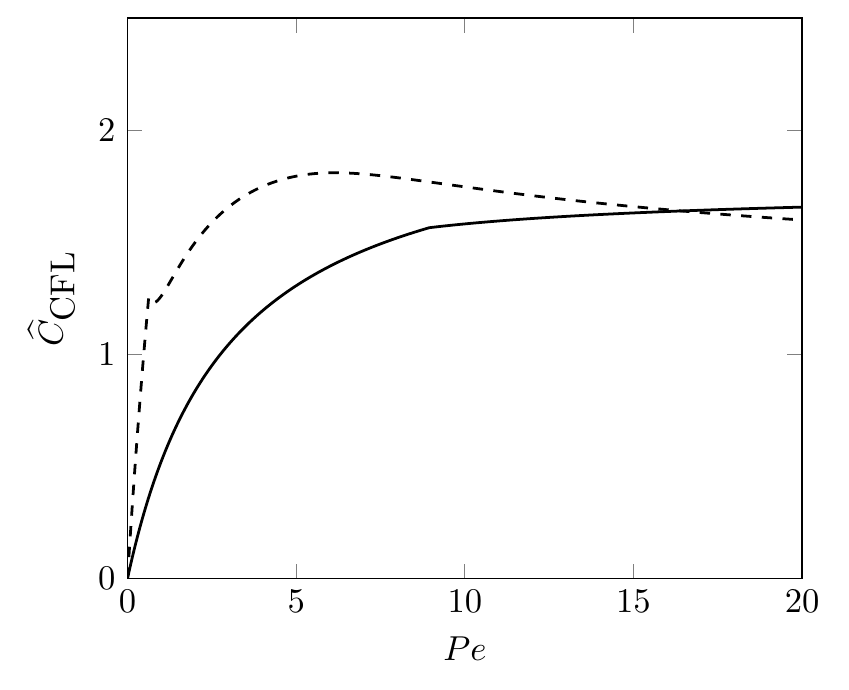}
\end{tabular}
\caption{$Pe\to \widehat C_{\text{CFL}}$ curves for the four
scenarios that we are considering: centered (left) and weak upwind
(right), with solid --- RK$_4$, dashed --- RK$_\text{D}$.}
\label{fig:PeVsCFL}
\end{figure}

The curves $Pe\to \widehat
C_{\text{CFL}}(Pe;\theta_3,\theta_4,w_3,w_4)$ are determined for the
four situations we want to highlight: RK$_4$ and RK$_\text{D}$ for
time; centered and weak upwind for space. We plot in
Figure~\ref{fig:PeVsCFL} the optimal CFL curves for the fourth-order
centered scheme and the third order weak upwind scheme with $Pe\in
[0.001,20]$ (note that the oscillations we observe for the centered
case are a consequence of the discrete consideration of the spectral
curve). We complement the figure with Table~\ref{tab:CFL_big_Pe}
indicating the optimal CFL values, $\widehat C_{\text{CFL}}$, for
$Pe=0$ (an expression depending on $\Delta x$) and large P\'eclet
values.

\begin{table}[!ht]\centering
\caption{Comparison of $\widehat C_\text{CFL}$ for RK$_{\text{D}}$
and RK$_4$.} \label{tab:CFL_big_Pe}
\begin{tabular}{@{}ccccc@{}}\toprule
$Pe$ & \multicolumn{2}{c}{RK$_{\text{D}}$} &
\multicolumn{2}{c}{RK$_4$}\\\cmidrule(l{3pt}r{3pt}){2-3}
\cmidrule(l{3pt}r{3pt}){4-5} & centered & weak upwind & centered &
weak upwind \\\midrule
0 & $1.77\Delta x^2$ & $2.36\Delta x^2$ & $0.53\Delta x^2$ & $0.70\Delta x^2$ \\
20 & 1.04 & 1.60 & 1.62 & 1.66 \\
200 & 0.45 & 1.34 & 2.04 & 1.74 \\
20000 & 0.09 & 1.25 & 2.06 & 1.75 \\
200000 & 0.04 & 1.24 & 2.06 & 1.75 \\
$+\infty$ & 0 & 1.24 & 2.55 & 1.75 \\ \bottomrule
\end{tabular}
\end{table}

From Figure~\ref{fig:PeVsCFL} and Table~\ref{tab:CFL_big_Pe} we draw
the following conclusions:
\begin{itemize}
\item \textbf{centered scheme} RK$_\text{D}$ provides the largest time-steps for $Pe<5$ while RK$_4$ is more efficient for $Pe>5$. In particular, the limit for $Pe\to +\infty$ of $\widehat C_\text{CFL}$ is 2.55 for RK$_4$ while $\widehat C_\text{CFL}$ converges to 0 with RK$_\text{D}$.
\item \textbf{weak upwind scheme} For $Pe<15$, $\widehat C_\text{CFL}$ value is larger with RK$_\text{D}$ scheme than RK$_4$. The par\-ti\-cu\-lar case $Pe=0$ shows that RK$_\text{D}$ case allows a time parameter about three times larger than the RK$_4$ case. For $Pe>15$, RK$_4$ method turns out to be more efficient.
\end{itemize}

\subsubsection{Convergence order and stability}

\begin{figure}[!ht]\centering
\begin{tabular}{@{}ccc@{}}\toprule
$I$ & centered & weak upwind \\\midrule 25 &
\includegraphics[align=c,height=4cm]{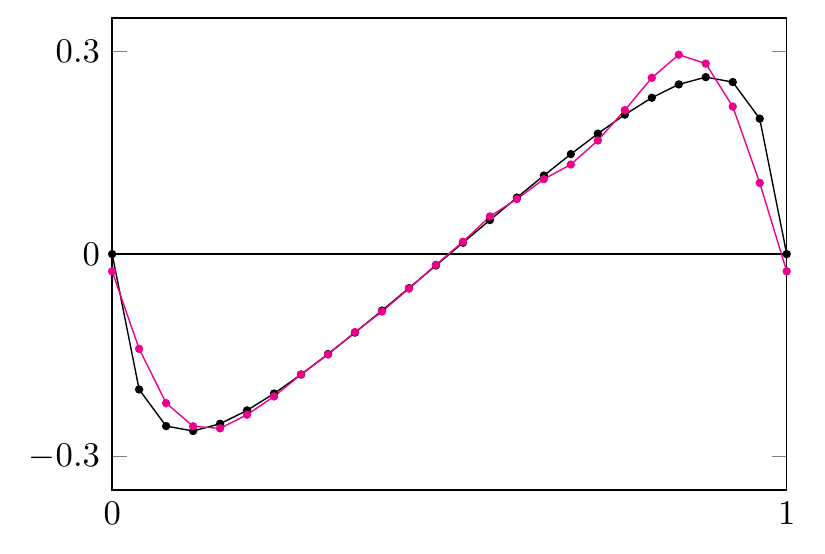} &
\includegraphics[align=c,height=4cm]{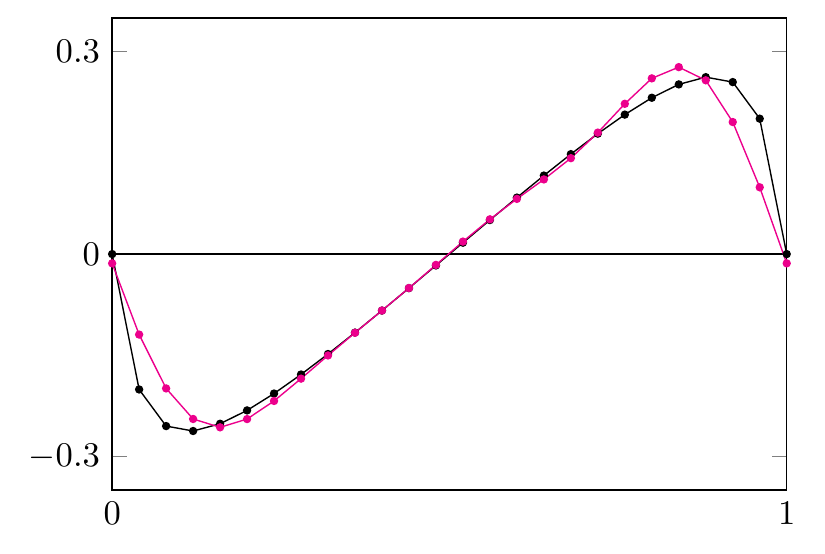}\\\midrule
50 &
\includegraphics[align=c,height=4cm]{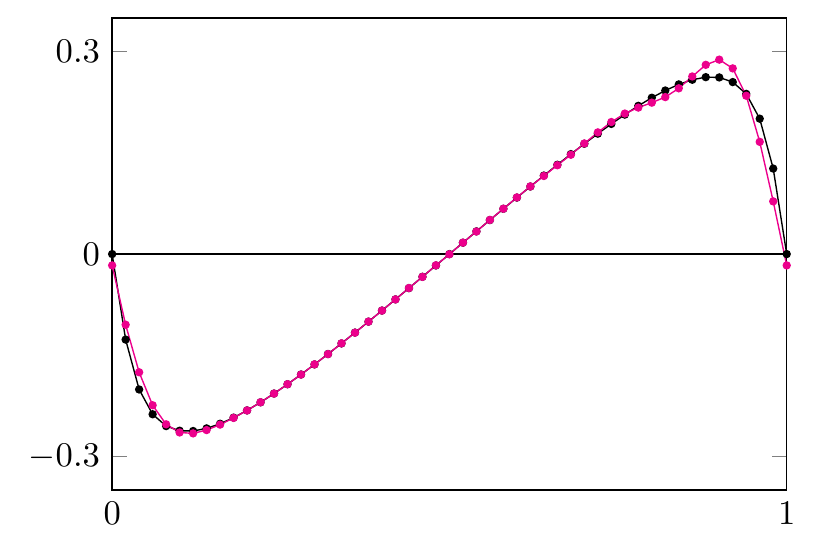}
&
\includegraphics[align=c,height=4cm]{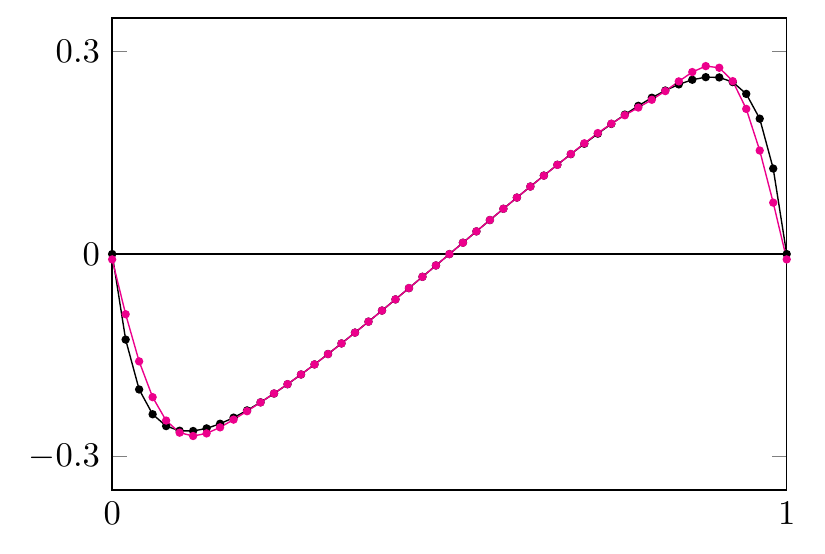}\\\midrule
100 &
\includegraphics[align=c,height=4cm]{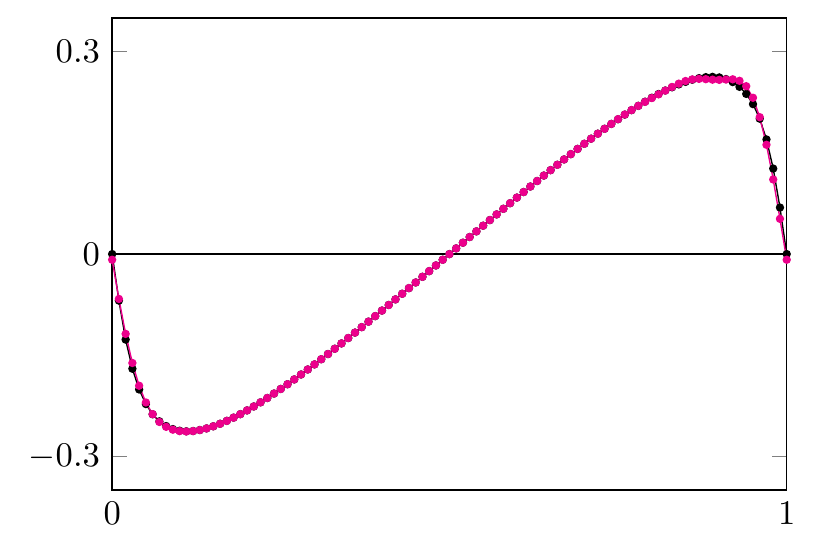}
&
\includegraphics[align=c,height=4cm]{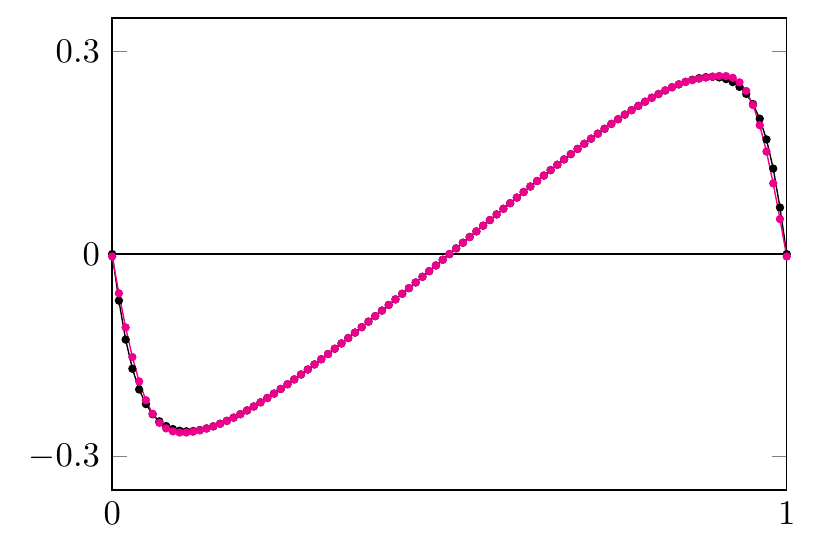}\\\bottomrule
\end{tabular}
\caption[]{Exact (\exactsol) and approximate (\numericalsol)
solutions considering RK$_4$ time scheme for benchmark~4.}
\label{fig:compare_RK4_centered_upwind_regular}
\end{figure}

To compare the convergence error between the centered and the weak
upwind discretisations, we consider the advection ($u=1$,
$\kappa=0$, $Pe=+\infty$) of the manufactured
solution~\eqref{eq:delta_solution_xt} with $\delta=0.1$,
corresponding to a smooth solution since the transition takes place
in more than 5 nodes. \vskip 1em
--- benchmark~4.
We carry out the RK$_4$ scheme in time with a time step $\Delta
t=0.8\Delta t_{\max}$ until $t_{\text{f}}=1$. We present in
Table~\ref{tab:max_time_step_centered} the time step, error
$E_\infty$
\[
E_{\infty}\equiv E_\infty(\Phi^N,I) =\displaystyle\max_{i=1}^{I}
|\phi_i^N-\phi(x_i,t_\text{f})|,
\]
and the respective convergence order $O_\infty$ between two
solutions/grids $(\Phi_k,I_k)$, for $k=1,2$ where $I_1<I_2$ as
\[
O_{\infty} \equiv O_{\infty} \left( (\Phi_1^N,I_1),(\Phi_2^N,I_2)
\right) = \frac{|\log
E_\infty(\Phi_1^N,I_1)/E_\infty(\Phi_2^N,I_2)|}{|\log I_1/I_2|},
\]
while we display in
Figure~\ref{fig:compare_RK4_centered_upwind_regular} the numerical
approximations for the centered and weak upwind schemes. We reach
the fourth-order in space for the centered method and the expected
third-order in space for the upwind case (the fourth-order method in
time turns to be insignificant for large values of $I$, since small
values of $\Delta t$ are needed due to the CFL condition).

\begin{table}[!ht]\centering
\caption{Time steps and errors using 80\% of the maximal time step
given by $\widehat C_\text{CFL}$ considering RK$_4$ time scheme  for
benchmark~4.} \label{tab:max_time_step_centered}
\begin{tabular}{@{}ccccccc@{}}\toprule
$I$ & \multicolumn{3}{@{}c@{}}{centered} &
\multicolumn{3}{@{}c@{}}{weak
upwind}\\\cmidrule(l{3pt}r{0pt}){2-4}\cmidrule(l{3pt}r{0pt}){5-7}
 & $\Delta t$ & $E_\infty$ & $O_\infty$ & $\Delta t$ & $E_\infty$ & $O_\infty$ \\\midrule
100  & 1.65E$-$2 & 1.68E$-$2 & --- & 1.40E$-$2 & 2.20E$-$2 & --- \\
200  & 8.25E$-$3 & 3.88E$-$3 & 2.1 & 6.98E$-$3 & 6.29E$-$3 & 1.8 \\
400  & 4.12E$-$3 & 4.88E$-$4 & 3.0 & 3.49E$-$3 & 1.20E$-$3 & 2.4 \\
800  & 2.06E$-$3 & 3.83E$-$5 & 3.7 & 1.75E$-$3 & 1.70E$-$4 & 2.8 \\
1600  & 1.03E$-$3 & 2.46E$-$6 & 4.4 & 8.73E$-$4 & 2.15E$-$5 & 3.0
\\\bottomrule
\end{tabular}
\end{table}

To check the stability condition, we have performed the computation
with the RK$_4$ time scheme with three different time steps: (1) the
time step $\Delta t_{\max}$ corresponding to the optimal CFL value
given by $\widehat C_\text{CFL}(Pe)$, (2) a $20\%$ smaller time
step, and (3) a $10\%$ larger time step. We present in
Table~\ref{tab:timestep_compare} the errors $E_\infty$ together with
the number of iterations $n_\text{TS}$ needed to reach the final
time $t_{\text{f}}=1$. For the latter case, stability is no longer
preserved and the error blows up. With the critical time step
$\Delta t_{\max}$, we manage to compute the solution until the final
time but with an error slightly larger than the one obtained with
the smaller $\Delta t$. Notice that, as expected, the number of
steps linearly increases with the number of nodes.

\begin{table}[!ht]\centering
\caption{Stability study considering RK$_\text{4}$ time scheme for
benchmark~4.} \label{tab:timestep_compare}
\begin{tabular}{@{}cccccccccc@{}}\toprule
$I$  & space scheme & $\widehat C_\text{CFL}$ & $\Delta t_{\max}$ &
\multicolumn{2}{c}{$\Delta t=\Delta t_{\max}$} &
\multicolumn{2}{c@{}}{$\Delta t=0.8\Delta t_{\max}$} &
\multicolumn{2}{c@{}}{$\Delta t=1.1\Delta
t_{\max}$}\\\cmidrule(l{3pt}r{3pt}){5-6}\cmidrule(l{3pt}r{3pt}){7-8}\cmidrule(l{3pt}r{0pt}){9-10}
 & & & & $E_{\infty}$ & $n_\text{TS}$ & $E_{\infty}$ & $n_\text{TS}$ & $E_{\infty}$ & $n_\text{TS}$\\\midrule
\multirow{2}{*}{25} & centered & 2.06 & 8.25E$-$2 & 1.13E$-$1 & 13 & 9.52E$-$2 & 16 & 1.01E$+$1 & 12\\
& weak upwind & 1.77 & 7.06E$-$2 & 1.12E$-$1 & 15 & 1.02E$-$1 & 18 &
4.62E$-$1 & 13\\\midrule
\multirow{2}{*}{50} & centered & 2.06 & 4.13E$-$2 & 5.44E$-$2 & 25 & 4.86E$-$2 & 31 & 2.35E$+$3 & 23\\
& weak upwind & 1.75 & 3.49E$-$2 & 5.46E$-$2 & 29 & 5.05E$-$2 & 36 &
4.02E$+$0 & 27\\\midrule
\multirow{2}{*}{100} & centered & 2.06 & 2.06E$-$2 & 2.23E$-$2 & 49 & 1.68E$-$2 & 61 & 1.44E$+$8 & 45\\
& weak upwind & 1.75 & 1.75E$-$2 & 2.45E$-$2 & 58 & 2.20E$-$2 & 72 &
4.79E$+$2 & 53\\\midrule
\multirow{2}{*}{200} & centered & 2.06 & 1.03E$-$2 & 5.87E$-$3 & 98 & 3.88E$-$3 & 122 & NaN & 89\\
& weak upwind & 1.75 & 8.73E$-$3 & 7.05E$-$3 & 115 & 6.29E$-$3 & 144
& 1.03E$+$7 & 105\\\bottomrule
\end{tabular}
\end{table}

\subsection{Hybrid scheme}
Hybrid time scheme consists in choosing between the RK$_4$ and the
RK$_\text{D}$ scheme in function of the P\'eclet number to provide
the largest time step. For space (and time) dependent velocity and
diffusion coefficient, the scheme then would be different from one
node to another. For example, for the centered scheme in space, the
left panel of Figure~\ref{fig:PeVsCFL} shows that the highest
$\widehat C_{CFL}(Pe)$ values are reached with the RK$_\text{D}$
scheme when $Pe<5$ while the $\widehat C_{CFL}(Pe)$ is larger with
RK$_4$ when $Pe>5$. We implement an hybrid scheme that switches from
one method to the other in function of the velocity and diffusion
coefficients to optimise the time step and reduce the computational
effort.

We consider the space-dependent parameter equation $\partial_t
\phi+\mathfrak E(x)[\phi]=f(x)$ with
\begin{equation}
\label{space_ConvDiff_operator} \mathfrak
E(x)[\phi]=-u(x)\phi'+\kappa(x)\phi''
\end{equation}
and $f$ a given source term. Since the cell P\'eclet number depends
on the position, we set $\displaystyle Pe_i=\frac{u_i\Delta
x}{\kappa_i}$ for node $i$ with $u_i=u(x_i)$ and
$\kappa_i=\kappa(x_i)$ for the velocity and diffusion, respectively.
Taking the centered scheme for the discretization in space of
operator~\eqref{space_ConvDiff_operator}, the hybrid scheme in time
is obtained with the following rule:
\begin{itemize}
\item if $Pe_i<5$, we use the RK$_\text{D}$ scheme for node $i$;
\item otherwise, we use the RK$_4$ scheme.
\end{itemize}
Note that the two time schemes are compatible since we have the same
sub-steps by construction.

The weak upwind scheme in space is also considered and the hybrid
time scheme derives from the right panel of
Figure~\ref{fig:PeVsCFL}, and is given by:
\begin{itemize}
\item if $Pe_i<15$, we use the RK$_\text{D}$ scheme for node $i$;
\item otherwise, we use the RK$_4$ scheme.
\end{itemize}

To update the solution from $t^n$ to $t^{n+1}=t^n+\Delta t$, we
compute the global time step $\Delta t$ in the following way: the
P\'eclet number $Pe_i$ and the corresponding optimal CFL number
$(\widehat C_{CFL})_i$ are computed leading to an optimal time step
$\Delta t_i$ that provides the stability. In order to guarantee the
global stability, we then choose
\begin{equation}
\label{eq:time_parameter} \Delta t=\min_{i=1,\ldots,I} \Delta t_i.
\end{equation}

To test the hybrid scheme, we consider the manufactured solution
$\phi(x,t;\omega)=\sin(2\pi \omega(x-u t))$, with $\omega=1$ where
the velocity $u=1$ is constant and the diffusion is given by
\[
\kappa(x)=a_0\exp(25(x-1/2)^2)+a_1, \quad a_0=10^{-4}, \quad
a_1=10^{-5},
\]
\vskip 1em
--- benchmark~5.
The source term is calculated to satisfy relation $\partial_t
\phi=\mathfrak E(x)[\phi]+f(x)$ from the manufactured solution.
Numerical simulations are carried out until the final time
$t_{\text{f}}=1$ with $I=100$. It is worth noting that at the same
time step or sub-step, we handle two different schemes in time
depending on the cell P\'eclet number, with different Butcher's
tableaux.

\begin{table}[!ht]\centering
\caption{Error and number of iterations for the full RK$_4$ and
hybrid scheme for benchmark~5.} \label{tab:hybrid_scheme}
\begin{tabular}{@{}ccccccc@{}}\toprule
scheme & \multicolumn{3}{@{}c@{}}{centered} &
\multicolumn{3}{@{}c@{}}{weak
upwind}\\\cmidrule(l{3pt}r{0pt}){2-4}\cmidrule(l{3pt}r{0pt}){5-7}
 & $\Delta t$ & $E_\infty$ & $n_\text{TS}$ & $\Delta t$ & $E_\infty$ & $n_\text{TS}$ \\\midrule
 hybrid (RK$_\text{4}$/RK$_\text{D}$) & 3.50E$-$3  & 5.56E$-$5 & 286  & 4.38E$-$3  & 1.46E$-$4  & 229\\
RK$_\text{4}$  & 1.01E$-$3 & 3.03E$-$6 & 993 & 1.26E$-$3 & 7.80E$-$5 & 792\\
\bottomrule
\end{tabular}
\end{table}
We compare the hybrid scheme with the full RK$_4$ time scheme (which
does not depend on the P\'eclet number) and report in
Table~\ref{tab:hybrid_scheme} the comparison between the two
methods. On the one hand, the full RK$_4$ scheme achieves the better
accuracy (errors cut almost by twenty for the centered case and
almost by two for the weak upwind) but the hybrid scheme provides
the largest time steps with a strong reduction of the number of
iterations (almost four times faster). We highlight that no
oscillations appear and stability is achieved for both schemes as
shown by the errors convergence rate.

\begin{figure}[!ht]\centering
\begin{tabular}{@{}cc@{}}
\includegraphics[align=c,height=6.5cm]{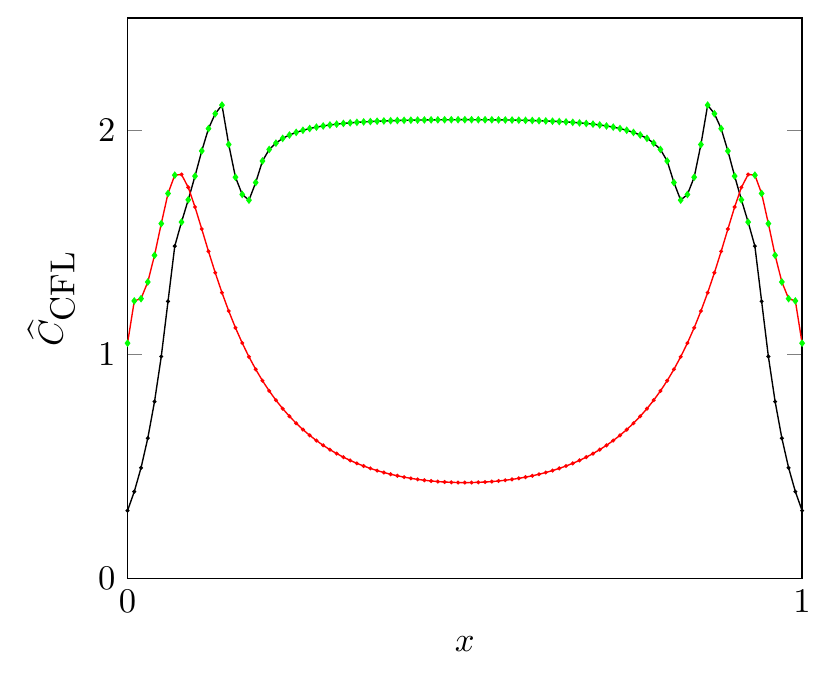} &
\includegraphics[align=c,height=6.5cm]{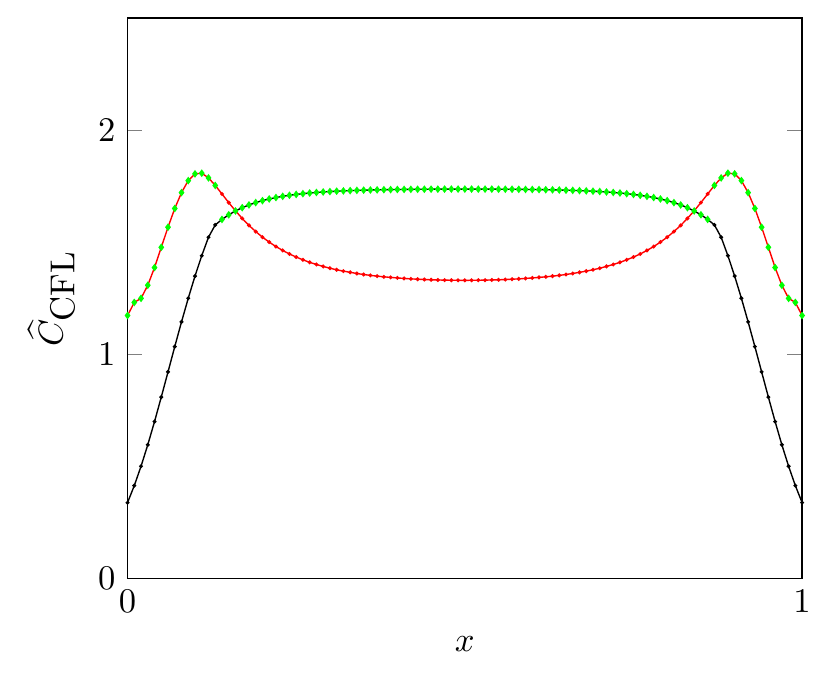}
\end{tabular}
\caption[]{$\widehat C_{CFL}$ plots for benchmark~5: centered (left)
and weak upwind (right).
 $\widehat C_{\text{CFL},\text{hybrid}}$ values (\hybridCFL);  $\widehat C_{\text{CFL},\text{RK}_4}$ curve (~\exactsol), and $\widehat C_{\text{CFL},\text{RK}_\text{D}}$ curve (\rkdCFL).}
\label{fig:CFL_hybrid}
\end{figure}

We display in Figure~\ref{fig:CFL_hybrid} the $\widehat
C_{\text{CFL},\text{RK}_4}$ calculated with the full RK$_4$ scheme
and the $\widehat C_{\text{CFL},\text{RK}_\text{D}}$ with the
RK$_\text{D}$ scheme while we highlight the $\widehat
C_{\text{CFL},\text{hybrid}}$ for the hybrid  scheme with the green
mark. We observe that, except for a small number of nodes, we have
the property
\[
\widehat C_{\text{CFL},\text{hybrid}}=\max\Big(\widehat
C_{\text{CFL},\text{RK}_4}, \widehat
C_{\text{CFL},\text{RK}_\text{D}}\Big)
\]
that confirms we have taken the best scheme in time for each node
providing the larger global time step~$\Delta t$.

\section{The \textit{a posteriori} method for optimal time step scheme}

At the numerical level, smooth solutions are functions where the
numerical derivative is bounded for $\Delta x$ small enough and
transitions between successive extremes are spread on, at least,
five or six nodes. In that case, even with a large P\'eclet number,
the centered scheme is stable and the time step is ruled by the
hybrid scheme condition. On the other hand, to handle sharp
gradients or even discontinuities, the weak upwind scheme turns out
to be the candidate method substituting the centered scheme,
possibly leading to a change of the scheme in time. At time $t^n$
and for each node $i$, the scheme in space has to be chosen with
respect to the local regularity for that particular node. Then the
scheme in time is chosen with respect to the scheme in space, the
local cell P\'eclet, and $\widehat C_{CFL}$  to provide the optimal
local time step $\Delta t^n_i$. Then, using the time parameter
$\Delta t^n$ given by relation (\ref{eq:time_parameter}), we update
the solution at time $t^{n+1}$.

In order to make the optimal choice for the space and time schemes
(node by node), we use the \textit{a posteriori} paradigm also
mentioned as MOOD method for Multi-dimensional Optimal Order
Detector \cite{Clain, Clain2}. In this context, the strategy is
based on the choice between a high accurate scheme in space (the
centered scheme) and a high stable scheme (the weak upwind scheme)
together with the associated optimal time scheme we studied in the
hybrid scheme section.

\subsection{Basics on \textit{a posteriori} strategy}

We present a brief description of the \textit{a posteriori} paradigm
and introduce the notations we need in this section. Assume that the
numerical solution $\Phi^n$ is known at time $t^n$.
\begin{enumerate}
\item We compute a candidate solution $\Phi^\star$ for time $t^{n+1}$ using the most accurate schemes in space and time, namely the centered scheme and the RK$_4$.
\item We check, node by node, the admissibility of the solution using detectors, \textit{i.e.} small routines that analyse specific aspects of the approximations such as extrema, oscillations, and physical property violations if any.
\item The nodes detected as non-admissible are computed again but with the weak upwind and time schemes using the hybrid scheme procedure.
\end{enumerate}
\begin{myremark}
The method is tagged \textit{a posteriori} since we analyse a
candidate solution $\Phi^\star$ after computing a time step. On the
contrary, MUSCL and WENO method are said \textit{a priori} since we
perform the limitation strategy based on solution $\Phi^n$,
\cite{MUSCL, WENO}.
\end{myremark}

\subsubsection{Node space and time scheme tables}

In practice, we introduce two tables \texttt{CSS} and \texttt{CTS}
for Cell Space Scheme and Cell Time Scheme, respectively, with the
following rules.
\begin{itemize}
\item We set $\texttt{CSS}[i]=0$ if we use the centered scheme otherwise $\texttt{CSS}[i]=1$ for the weak upwind scheme.
\item We set $\texttt{CTS}[i]=0$ if we use the RK$_4$ scheme otherwise $\texttt{CTS}[i]=1$ for the RK$_\text{D}$ scheme.
\end{itemize}
Given the numerical solution $\Phi^{\star,0}=\Phi^{n}$ and tables
\texttt{CSS} and \texttt{CTS}, we rewrite the RK scheme taking the
two tables into account. For each node $i$ the terms
$\Phi_i^{\star,j}$ and $\mathcal K_i^{\star,j}$  are computed with
\begin{align*}
& \Phi_i^{\star,j}=\Phi_i^{\star,0}+\Delta t\sum_{\ell=1}^{s}a_{j\ell}(\texttt{CTS}[i]) \mathcal K_i^{\star,\ell},\\
& \Phi_i^{\star}=\Phi_i^{\star,0}+\Delta t \sum_{j=1}^{s}b_j(\texttt{CTS}[i]) \mathcal K_i^{\star,j},\\
& \mathcal K_i^{\star,j}=E_i(\texttt{CSS}[i])
\Phi^{\star,j}+F(x_i,t^{\star,j}), \quad j=1,\ldots,s.
\end{align*}
Expressions $a_{j\ell}(\texttt{CTS}[i])$ and $b_j(\texttt{CTS}[i])$
indicate that we use RK$_4$ scheme if $\texttt{CTS}[i]=0$ or
RK$_\text{D}$ scheme if $\texttt{CTS}[i]=1$. On the other hand, the
term $E_i(\texttt{CSS}[i])$ states that we use the centered scheme
for node $i$ if $\texttt{CSS}[i]=0$ or the weak upwind scheme if
$\texttt{CSS}[i]=1$.

\subsubsection{Detectors}

Detectors  are small routines to check a specific property of the
candidate solution. We assemble the detectors in a chain of
operations which basically indicate if a node would be cured (change
the scheme in space) or not. We list hereafter the detectors we use
in the present document and refer to \cite{Clain3} for a detailed
presentation of the most useful detectors.

$\bullet$ \textbf{ED.} The Extrema Detector intends to localise the
extrema of the numerical function by checking
\[
s_i=\text{sign}\Big((\phi_i^{\star}-\phi_{i+1}^{\star})(\phi_i^{\star}-\phi_{i-1}^{\star})\Big).
\]
If $s_i>0$ we have an extremum and the detector is activated and
returns \texttt{true} otherwise it corresponds to a monotone
situation and the detector returns \texttt{false}.

$\bullet$ \textbf{SCD.} The Small Curvature Detector helps to select
the very small oscillations that we consider innocuous from the
stability point of view. We calculate the variation quantity
\[
v_i=\max \left (\frac{|\phi_i^{\star}-\phi_{i+1}^{\star}|}{\Delta
x},\frac{|\phi_i^{\star}-\phi_{i-1}^{\star}|}{\Delta x}\right )
\]
and, for a user parameter $\theta >0$, the detector returns
\texttt{true} if $v_i < \theta \Delta x$ and \texttt{false}
otherwise since a small value of $v_i$ indicates oscillations with
too low magnitude to be considered as an issue.

$\bullet$ \textbf{LOD.} Local Oscillation detection aims to detect
variations deriving from oscillations. Indeed, a local oscillation
is characterised by a variation of the curvature sign. To this end,
we compute the second derivative
\[
\chi_i=\frac{\phi_{i+1}^{\star}-2\phi_i^{\star}+\phi_{i-1}^{\star}}{\Delta
x}.
\]
Then the detector is deactivated (return \texttt{false}) if
$\chi_{i-1}$, $\chi_i$, and $\chi_{i+1}$ have the same sign and is
activated (return \texttt{true}) if one of the curvatures has a
different sign of the two others.

$\bullet$ \textbf{SD.} Smooth Detector consists in assessing the
local numerical smoothness to determine if an extremum is physical
or if it corresponds to a numerical artefact. To evaluate the
smoothness of the solution we use, once again, the curvature and
define the minimum and maximum absolute curvature
\[
\chi_{i,m}=\min(\chi_{i-1},\chi_{i},\chi_{i+1}),\qquad
\chi_{i,M}=\max(\chi_{i-1},\chi_{i},\chi_{i+1}).
\]
For a given user parameter $\theta\in [0,1]$, the solution is not
smooth (detector activated \texttt{true}) if $\chi_{i,m}<\theta
\chi_{i,M}$ since we detect large variation of curvature on the
three consecutive points. Otherwise, the detector returns
\texttt{false} that indicates the solution is considered smooth
enough at the numerical level.

\subsubsection{Detector chain and \textit{a posteriori} cure}

The detectors being defined, we assemble it into an ordered chain
that enables to decide whether a node would be corrected or not.
\begin{figure}\centering
\includegraphics[width=0.8\textwidth]{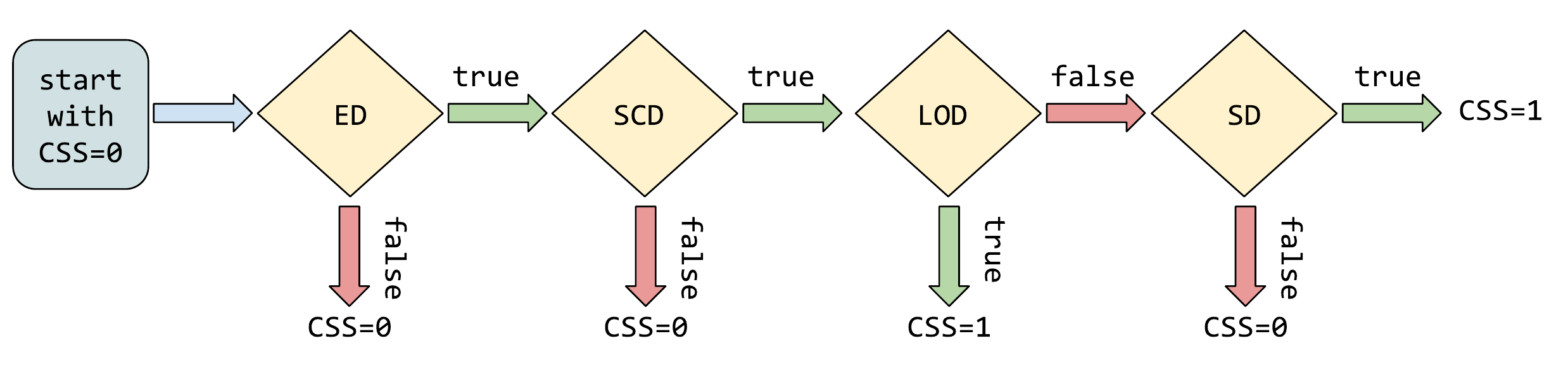}
\caption{The chain detector.} \label{fig:chain_detector}
\end{figure}
The detector chain is carried out for nodes with
$\texttt{CSS[i]=0}$, that is the centered scheme. We preserve the
high accuracy scheme if there is no extremum (ED \texttt{false}) or
generated by too small variations (SCD \texttt{false}). When a
potential oscillation is detected (ED \texttt{true}), we try to
relax the scheme (and preserve the accuracy) with other detectors
that assess if the extremum is a real physical one. The LOD is
activated if we observe the sign change of the curvature and then we
set $\texttt{CSS}[i]=1$ to indicate that the node has to be cured.
In the same way, a node with too large variations of the curvature
is considered as a problematic node and we set $\texttt{CSS}[i]=1$
if SD is \texttt{true}. At the end of the chain, we have a new
\texttt{CSS} map that indicates the points we have to compute again
with the weak upwind scheme. The time scheme is determined with the
hybrid scheme strategy and we flag the $\texttt{CTS}[i]$
accordingly.

From the \texttt{CSS} map and the cell P\'eclet number, we modified
the \texttt{CTS} map following the rule given by the hybrid scheme.
At the end of the day, we get a new $\Delta t^n_i$ for each node $i$
and use the minimum time step for $\Delta t^n$ following
(\ref{eq:time_parameter}).

\begin{myremark}
Computational resources are reduced by only computing again the node
$i$ that have been cured together with the neighbour nodes that may
be affected by the values of $\Phi_i^{\star}$ during the four-stage
Runge-Kutta procedure. In practice, very few nodes are modified by
the \textit{a posteriori} correction (less than 5\%, see
\cite{Clain2}) and the additional cost is of the same order.
\end{myremark}

\subsection{Numerical tests}
We present two numerical tests to examine the efficiency of the
\textit{a posteriori} method. The first sanity check consists in
carrying out a simulation with a regular solution. Indeed, for
smooth approximations, the limiter strategy has to preserve the
higher accuracy and the chain detector has to return
$\texttt{CSS}[i]=0$ for all nodes $i$. \vskip 1em
--- benchmark~6.
We consider the manufactured solution~\eqref{eq:delta_solution} with
$\delta=0.15$. We take constant physical parameters $u=1$,
$\kappa=2.7778$E$-03$ and a grid of $I=60$ nodes to obtain the cell
P\'eclet number $Pe=6$. The simulation is carried out  until the
final time $t_{\text{f}}=0.5$, corresponding to half a revolution

\begin{figure}[ht!]\centering
\begin{tabular}{@{}cc@{}}
\includegraphics[align=c,height=4cm]{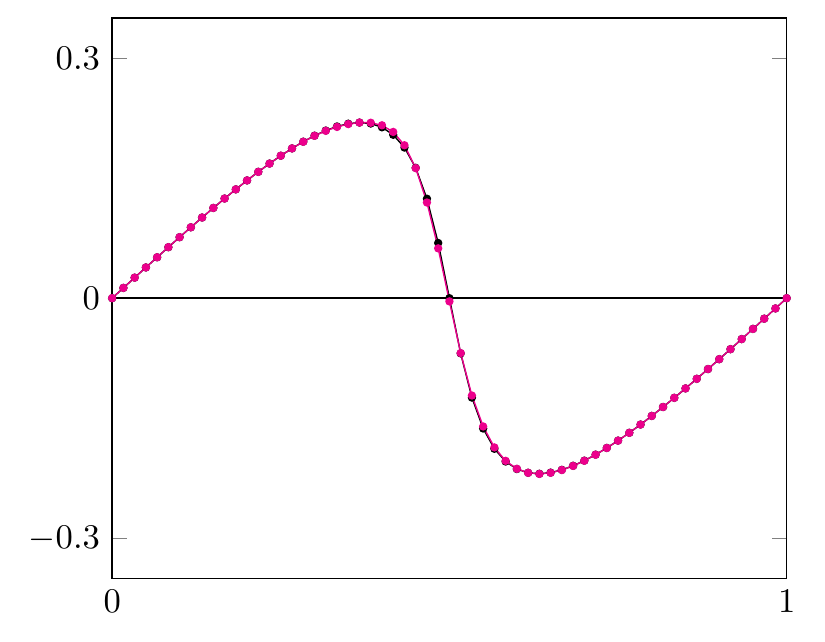} &
\includegraphics[align=c,height=4cm]{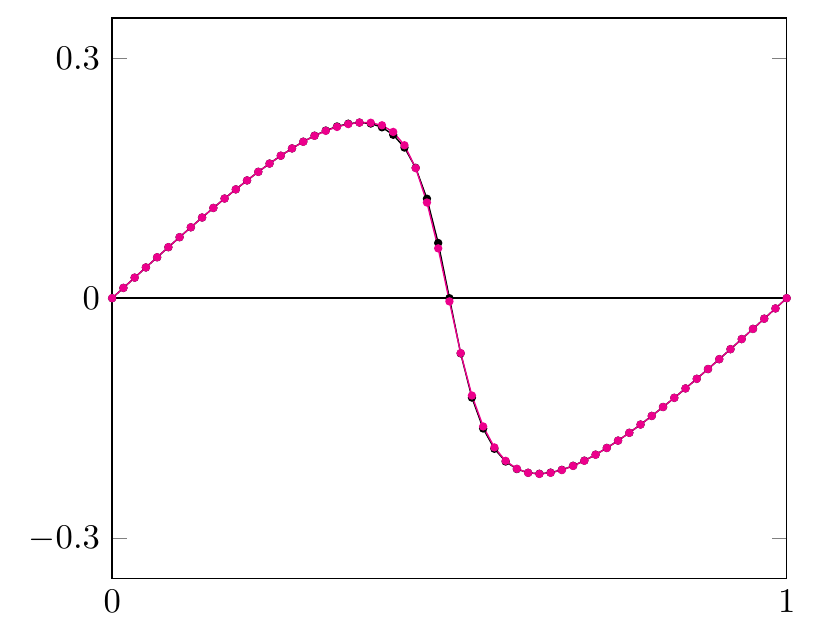} \end{tabular}
\caption[]{Exact (\exactsol) and approximate (\numericalsol)
solutions for $t_\text{f}=0.5$ the full RK$_4$+centered scheme
(left) and the \textit{a posteriori} strategy (right): the smooth
solution case --- benchmark~6.}
\label{fig:aposteriori_regular_solution}
\end{figure}

We display in Figure~\ref{fig:aposteriori_regular_solution}-left the
numerical approximation with the most accurate scheme
(RK$_4$+centered) while we reproduce on the right side the
approximation computed with the \textit{a posteriori} strategy. We
have checked that the table \texttt{CSS} has never been altered
during the simulation, that is, computation have been achieved with
the centered scheme and the RK$_4$ scheme in time due to the
P\'eclet number. \vskip 1em
--- benchmark~7.
The last benchmark deals with a rough function (regarded to the
characteristic mesh size) using the manufactured
solution~\eqref{eq:delta_solution} but with $\delta=0.015$
corresponding to a steep variation we assimilate as a shock regarded
to the small number of nodes $I=60$. We take constant physical
parameters $u=1$ and $\kappa=5.5556$E$-03$ to obtain the cell
P\'eclet number $Pe=3$. The simulation is again carried out  until
the final time $t_{\text{f}}=0.5$, corresponding to half a
revolution

\begin{figure}[h!]\centering
\begin{tabular}{@{}cc@{}}
\includegraphics[align=c,height=4cm]{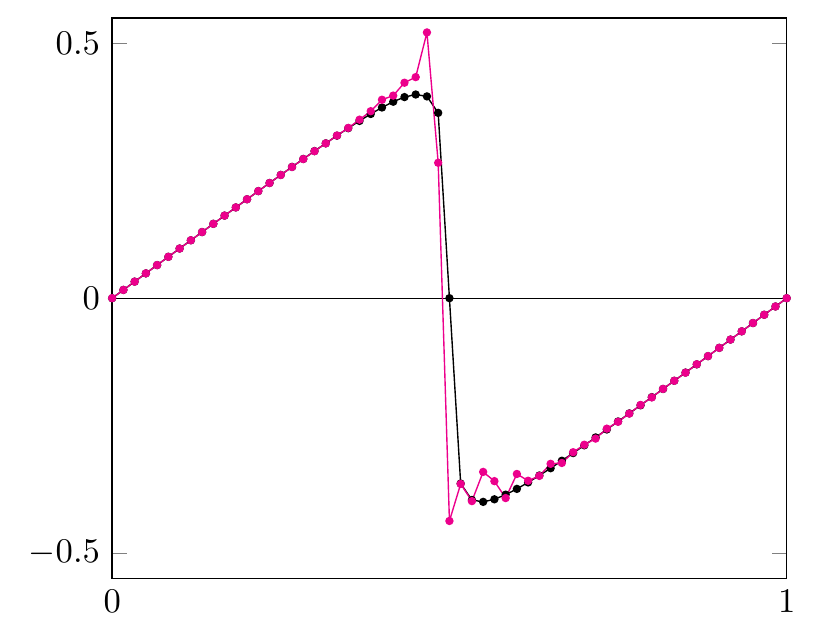} &
\includegraphics[align=c,height=4cm]{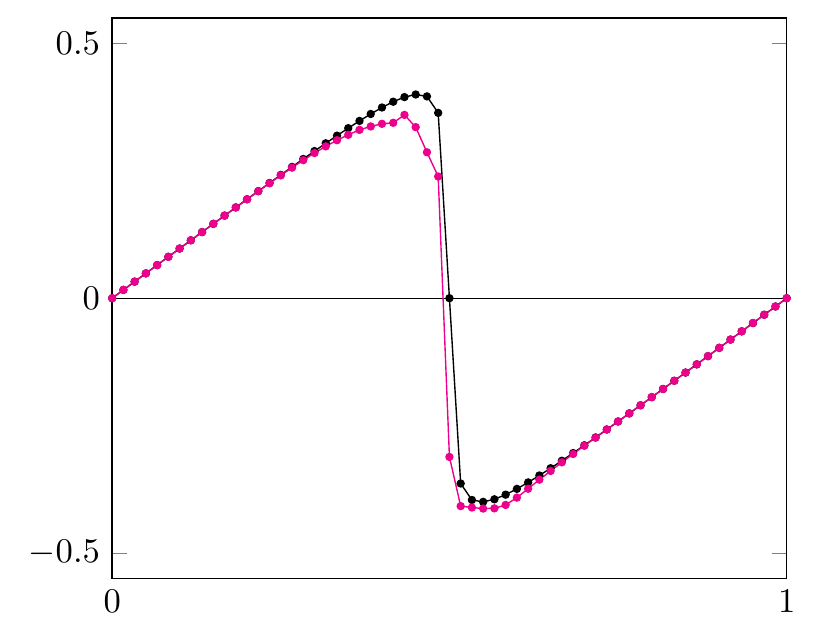} \end{tabular}
\caption[]{Exact (\exactsol) and approximate (\numericalsol)
solutions for $t_\text{f}=0.5$ the full RK$_4$+centered scheme
(left) and the \textit{a posteriori} strategy (right): the rough
solution case --- benchmark~7.}
\label{fig:aposteriori_rough_solution}
\end{figure}

We display in Figure~\ref{fig:aposteriori_rough_solution} the
solution obtained with the ``unlimited'' RK$_4$+centered scheme
(left) and the \textit{a posteriori} method (right). Clearly, the
steep gradient provokes oscillations when employing the low
diffusive centered scheme while the introduction of the weak upwind
scheme in some nodes (indicated with the red X on the figure)
manages to stabilise the solution and strongly reduces the over- and
under-shooting. Moreover the scheme in time on the node $i$ such
that $\texttt{CSS}[i]=1$ (weak upwind) switch to the RK$_\text{D}$
once $\texttt{CTS}[i]=1$ while the nodes where we maintain the
original centered scheme $\texttt{CSS}[i]=0$ (centered upwind) use
the RK$_4$ once $\texttt{CTS}[i]=0$. Notice that the number of nodes
that have been cured is almost 3 or 4, \textit{i.e.}, less than 8\%
for a 60-nodes grid.

\section{Conclusions and further work}
\label{sec:conclusions_and_further_work} We have developed a
strategy to analyse and optimise the stability based, on the one
hand, on the two-parameter family of continuous spectral curves that
characterise the space discretization and, on the other hand, a
two-parameter family of Runge-Kutta stability regions. Optimisation
results from the inclusion of the spectral curves into the stability
region with the help of the additional CFL parameter. We have
detailed the procedure with the five-point finite difference method
context but extension to other methods such as finite volume or
finite elements methods could be considered. A hybrid time scheme as
a function of the P\'eclet number have been proposed and analysed
with the objective of providing the largest time step while
preserving the stability. At last, we have presented an adaptation
of the \textit{a posteriori} strategy to handle the schemes in space
and time to preserve both the accuracy and the stability, even for
rough solutions, while we optimise the time step to reduce the
computational effort.

\section*{Acknowledgements}

G.J. Machado and S. Clain acknowledge the financial support by FEDER
-- Fundo Europeu de Desenvolvimento Regional, through COMPETE 2020
-- Programa Operational Fatores de Competitividade, and the National
Funds through FCT -- Funda\c c\~ao para a Ci\^encia e a Tecnologia,
project no. UID/FIS/04650/2019.

M.T. Malheiro acknowledge the financial support by Portuguese Funds
through FCT (Funda\c c\~ao para a Ci\^encia e a Tecnologia) within
the Projects UIDB/00013/2020 and UIDP/00013/2020 of CMAT-UM.

M.T. Malheiro, G.J. Machado, and S. Clain acknowledge the financial
support by FEDER -- Fundo Europeu de Desenvolvimento Regional,
through COMPETE 2020 -- Programa Operacional Fatores de
Competitividade, and the National Funds through FCT -- Funda\c c\~ao
para a Ci\^encia e a Tecnologia, project no.
POCI-01-0145-FEDER-028118.

\bibliographystyle{spmpsci}      % mathematics and physical sciences

\end{document}